\documentclass{article}
\usepackage{amsfonts}
\usepackage{amssymb}
\usepackage{graphicx}
\usepackage{amsmath}

\input{tcilatex}
\begin{document}

\begin{center}
\textbf{Singular perturbation for abstract elliptic equations and application%
}

\textbf{Veli Shakhmurov}\ \ 

Okan University, Department of Mechanical engineering, Akfirat, Tuzla 34959
Istanbul, Turkey, E-mail: veli.sahmurov@okan.edu.tr
\end{center}

\begin{quote}
\ \ \ \ \ \ \ \ \ \ \ \ \ \ 

\ \ \ \ \ \ \ \ \ \ \ \ \ \ \ \ \ \ \ \ \ \ \ \ \ \ \ \ \ \ \ \ \ \ \ \ \ \
\ \ \ \ \ \ \ \ \ \ \ \ \ \ \ 
\end{quote}

\begin{center}
\textbf{ABSTRACT}
\end{center}

\begin{quote}
\ \ \ \ \ \ \ \ \ \ \ \ \ \ \ \ \ 
\end{quote}

\ \ \ Boundary value problem for complete second order elliptic equation is
considered in Banach space. The equation and boundary conditions involve a
small and spectral parameter. The uniform $L_{p}-$regularity properties with
respect to space variable and parameters are established. Here, the explicit
formula for the solution is given and behavior of solution is derived when
the small parameter approaches zero. It used to obtain singular perturbation
result for abstract elliptic equation


\begin{center}
\bigskip\ \ \textbf{AMS: 34G10, 35J25, 35J70}

\textbf{1. Introduction, notations and background }
\end{center}

It is well known that differential equations with small parameter play
important role in modeling of physical processes. Differential-operator
equations (DOEs) with parameter have also significant applications in
nonlinear analysis. DOEs are studied in $\left[ 1,\text{ }2\right] $, $\left[
4-7\right] ,$ $\left[ 9-14\right] ,$ $\left[ 16-24\right] $ and the
references therein. Main aim of this paper is to show the uniform
separability properties of boundary value problems (BVPs) for elliptic DOE
with parameters\ \ 
\begin{equation}
\ -\varepsilon u^{\left( 2\right) }\left( t,\varepsilon \right) +Au\left(
t,\varepsilon \right) +Bu^{\left( 1\right) }\left( t,\varepsilon \right)
+\lambda u\left( t,\varepsilon \right) =f\left( t\right) ,  \tag{1.1}
\end{equation}%
where $A$, $B$ are linear operators in a Banach space $E,$ $\varepsilon $ is
a small and $\lambda $ is a complex parameter. Particularly, the sharp
coercive $L_{p}$ estimates for solution of $\left( 1.1\right) $ are obtained
uniformly with respect to small and spectral parameter. Finally, these
results are used in the singular perturbation problem, i.e. to study the
behavior of solution $u\left( t,\varepsilon \right) $ of ($1.1)$ and
convergence of $u\left( t,\varepsilon \right) $ as $\varepsilon \rightarrow
0 $ to the corresponding solution of the Cauchy problem for abstract
parabolic equation

\begin{equation}
\ Bu^{\left( 1\right) }\left( t\right) +Au\left( t\right) =f\left( t\right) ,
\tag{1.2}
\end{equation}%
\begin{equation*}
u\left( 0\right) =u_{0}.
\end{equation*}

The treatment of the singular perturbation problem for parabolic equation is
due to Fattorini $\left[ \text{7, Ch.VI}\right] $] (see also the references
therein). The singular perturbation problem for abstract hyperbolic equation 
\begin{equation}
\ \varepsilon u^{\left( 2\right) }\left( x,\varepsilon \right) +Au\left(
x,\varepsilon \right) =f\left( x,\varepsilon \right) ,  \tag{1.3}
\end{equation}%
was first considered by Kisynski $\left[ 12\right] $ in the case where $A$
is a self adjoint, positive definite operator on a Hilbert space. Latter,
Sova $\left[ 15\right] $ study the problem under the assumptions that $A$ is
the generator of a strongly continuous cosine function.

Then in $\left[ 6\right] $ the same problem considered for the complete
hyperbolic equation 
\begin{equation*}
\ \varepsilon u^{\left( 2\right) }\left( x,\varepsilon \right) +Au\left(
x,\varepsilon \right) +Bu^{\left( 1\right) }\left( x,\varepsilon \right) =0.
\end{equation*}

In contrast to these results, in this paper the singular perturbation
elliptic problem $\left( 1.1\right) $\ is considered and we show that the
solution $u\left( x,\varepsilon \right) $ of the equation $\left( 1.1\right) 
$ converge in $L_{p}\left( 0,1;E\right) $ as $\varepsilon \rightarrow 0$ to
the corresponding solution of the equation $\left( 1.2\right) $ uniformly
with respect to spectral parameter $\lambda .$ Moreover, the solution $%
u\left( \varepsilon ,x\right) $ of the elliptic BVP $\left( 1.1\right) $
converge in $E$ as $\varepsilon \rightarrow 0$ to the corresponding solution
of the Cauchy problem $\left( 1.2\right) $ uniformly with respect to
spectral parameter $\lambda .$ This result allow to investigate the spectral
properties of the parameter dependent elliptic BVP $\left( 1.1\right) .$
Since the Banach space $E$ is arbitrary and $A$ is a possible linear
operator, by chousing the spaces $E$ and operators $A$ we can obtained
different results about singular perturbation properties numerous classes of
elliptic, quasielliptic equations and its system which occur in a wide
variety of physical systems. Let we choose $E=L_{2}\left( 0,1\right) $ in $%
\left( 1.1\right) $ and $A$ to be differential operator with generalized
Wentzell-Robin boundary condition defined by 
\begin{equation*}
D\left( A\right) =\left\{ u\in W_{p_{1}}^{2}\left( 0,1\right) ,\text{ }%
Au\left( j\right) =0\text{, }j=0,1\right\} ,\text{ }
\end{equation*}%
\begin{equation*}
\text{ }Au=au^{\left( 2\right) }+bu^{\left( 1\right) }
\end{equation*}%
where $a$ is positive and $b$ is a real-valued functions. Assume $B$ is a
integral operator defined by 
\begin{equation*}
Bu=\dint\limits_{0}^{1}K\left( y,\tau \right) u\left( y,\tau \right) d\tau ,
\end{equation*}%
here, $K=K\left( y,\tau \right) $ is complex valued bounded function.

Then, we get the $L_{\mathbf{p}}\left( \Omega \right) -$separability and
singular perturbation properties of the Wentzell-Robin type BVP for elliptic
equation with integral term

\begin{equation}
-\left[ \varepsilon \frac{\partial ^{2}u}{\partial t^{2}}+a\frac{\partial
^{2}u}{\partial y^{2}}+b\frac{\partial u}{\partial y}\right]
+\dint\limits_{0}^{1}K\left( y,\tau \right) \frac{\partial }{\partial t}%
u\left( t,y,\tau \right) d\tau +\lambda u\left( t,y\right) =f\left(
t,y\right) \text{, }  \tag{1.4}
\end{equation}%
\begin{equation*}
\sum\limits_{i=0}^{m_{1}}\varepsilon ^{\frac{i}{2}}\alpha _{i}u^{\left(
i\right) }\left( 0,y,\varepsilon \right) =f_{1},\text{ }\sum%
\limits_{i=0}^{m_{2}}\varepsilon ^{\frac{i}{2}}\beta _{\imath }u^{\left(
i\right) }\left( T,y,\varepsilon \right) =f_{2}\text{ for a.e. }y\in \left(
0,1\right) ,
\end{equation*}%
\ \ \ 
\begin{equation}
a\left( j\right) u_{yy}\left( t,j,\varepsilon \right) +b\left( j\right)
u_{y}\left( t,j,\varepsilon \right) =0\text{, }j=0,1\text{, for a.e. }t\in
\left( 0,T\right) ,  \tag{1.5}
\end{equation}%
where $m_{k}\in \left\{ 0,1\right\} ,$ $\alpha _{i},$ $\beta _{i}$ are
complex numbers, $\varepsilon $ is a positive, $\lambda $ is a complex
parameter, $L_{\mathbf{p}}\left( \Omega \right) ,$ $\mathbf{p=}\left( p%
\mathbf{,}2\right) $ denotes mixed Lebesque space\ and $\Omega =\left(
0,T\right) \times \left( 0,1\right) $.

Note that, the regularity properties of Wentzell-Robin type BVP for elliptic
equations were studied e.g. in $\left[ \text{8}\right] $ and the references
therein.

We start by giving the notation and definitions to be used in this paper.

Let $E$ be a Banach space and $L_{p}\left( \Omega ;E\right) $ denotes the
space of strongly measurable $E$-valued functions that are defined on the
measurable subset $\Omega \subset R^{n}$ with the norm

\begin{equation*}
\left\Vert f\right\Vert _{L_{p}}=\left\Vert f\right\Vert _{L_{p}\left(
\Omega ;E\right) }=\left( \int\limits_{\Omega }\left\Vert f\left( x\right)
\right\Vert _{E}^{p}dx\right) ^{\frac{1}{p}}\text{, }1\leq p<\infty \ .
\end{equation*}

\ The Banach space $E$ is called $UMD$-space (see e.g. $\left[ 3\right] $)\
if the Hilbert operator 
\begin{equation*}
\left( Hf\right) \left( x\right) =\lim\limits_{\varepsilon \rightarrow
0}\int\limits_{\left\vert x-y\right\vert >\varepsilon }\frac{f\left(
y\right) }{x-y}dy
\end{equation*}%
is bounded in $L_{p}\left( R;E\right) $ for $p\in \left( 1,\infty \right) $. 
$UMD$ spaces include e.g. $L_{p}$, $l_{p}$ spaces and Lorentz spaces $L_{pq}$
for $p,$ $q\in \left( 1,\infty \right) $ and Morrey spaces (see e.g.$\left[
15\right] $ ).

Let $\mathbb{C}$ be the set of the complex numbers and\ 
\begin{equation*}
S_{\varphi }=\left\{ \lambda ;\text{ \ }\lambda \in \mathbb{C}\text{, }%
\left\vert \arg \lambda \right\vert \leq \varphi \right\} \cup \left\{
0\right\} ,\text{ }0\leq \varphi <\pi .
\end{equation*}

Let $B\left( E\right) $ denote the space of all bounded linear operators in $%
E$ and $R\left( \lambda ,A\right) $ denotes the resolvent of operator $A.$

A linear operator\ $A$ is said to be $\varphi $-positive in a Banach\ space $%
E$ with bound $M>0$ if $D\left( A\right) $ is dense on $E$ and 
\begin{equation*}
\ \left\Vert R\left( -\lambda ,A\right) \right\Vert _{B\left( E\right) }\leq
M\left( 1+\left\vert \lambda \right\vert \right) ^{-1}
\end{equation*}%
for any $\lambda \in S_{\varphi },$ $0\leq \varphi <\pi .$ Sometimes $%
A+\lambda I$\ will be denoted by $A+\lambda $ or $A_{\lambda },$ where $I$
denotes an identity operator in $E.$ It is known $\left[ \text{22, \S 1.15.1}%
\right] $ that there exist the fractional powers\ $A^{\theta }$ of a
positive operator $A.$ Let $E\left( A^{\theta }\right) $ denote the space $%
D\left( A^{\theta }\right) $ with norm 
\begin{equation*}
\left\Vert u\right\Vert _{E\left( A^{\theta }\right) }=\left( \left\Vert
u\right\Vert ^{p}+\left\Vert A^{\theta }u\right\Vert ^{p}\right) ^{\frac{1}{p%
}},1\leq p<\infty ,\text{ }0<\theta <\infty .
\end{equation*}

\ Let $E_{1}$ and $E_{2}$ be two Banach spaces. $\left( E_{1},E_{2}\right)
_{\theta ,p}$ for $0<\theta <1,1\leq p\leq \infty $ denotes the
interpolation spaces obtained from $\left\{ E_{1},E_{2}\right\} $ by the $K$%
-method $\left[ \text{22, \S 1.3.2}\right] $.\ 

$S\left( R^{n};E\right) $ is the Schwartz class, i.e. the space of all $E$%
-valued rapidly decreasing smooth functions on $R^{n}$ and $F$ denotes the
Fourier transformation. If the map $u\rightarrow \Lambda u=F^{-1}\Psi \left(
\xi \right) Fu,$ $u\in S\left( R^{n};E_{1}\right) $ is well defined and
extends to a bounded linear operator 
\begin{equation*}
\Lambda :\ L_{p}\left( R^{n};E_{1}\right) \rightarrow \ L_{p}\left(
R^{n};E_{2}\right)
\end{equation*}%
then a function $\Psi \in C\left( R^{n};B\left( E_{1},E_{2}\right) \right) $
is called a Fourier multiplier from $L_{p}\left( R^{n};E_{1}\right) $\ to $%
L_{p}\left( R^{n};E_{2}\right) .$

The set of all multipliers from\ $L_{p}\left( R^{n};E_{1}\right) $ to\ $%
L_{p}\left( R^{n};E_{2}\right) $\ will be denoted by $M_{p}^{p}\left(
E_{1},E_{2}\right) .$ For $E_{1}=E_{2}=E$\ it denotes by $M_{p}^{p}\left(
E\right) .$ Most important facts on Fourier multipliers and some related
reference can be found e.g. in $\left[ \text{22, \S 2.2.4}\right] $ and $%
\left[ \text{5, 23}\right] $.

Let 
\begin{equation*}
\Phi _{h}=\left\{ \Psi _{h}\in M_{p}^{p}\left( E_{1},E_{2}\right) ,\text{ }%
h\in Q\right\}
\end{equation*}%
be a collection of multipliers in $M_{p}^{p}\left( E_{1},E_{2}\right) $
dependent on the parameter $h.$ We say that $W_{h}$ is a uniform collection
of multipliers if there exists a positive constant $M$ independent on $h\in
Q $ such that

\begin{equation*}
\left\| F^{-1}\Psi _{h}Fu\right\| _{L_{p}\left( R^{n};E_{2}\right) }\leq
M\left\| u\right\| _{L_{p}\left( R^{n};E_{1}\right) }\ \ \ \ \ \ 
\end{equation*}
for all $h\in Q$ and $u\in S\left( R^{n};E_{1}\right) .$

Let $\mathbb{N}$, $\mathbb{R}$ denote the sets of natural and real numbers,
respectively. A set $G\subset B\left( E_{1},E_{2}\right) $ is called $R$%
-bounded (see e.g. $\left[ \text{5, 23}\right] $) if there is a positive
constant $C$ such that for all $T_{1},T_{2},...,T_{m}\in G$ and $%
u_{1,}u_{2},...,u_{m}\in E_{1},$ $m\in \mathbb{N}$ 
\begin{equation*}
\int\limits_{\Omega }\left\Vert \sum\limits_{j=1}^{m}r_{j}\left( y\right)
T_{j}u_{j}\right\Vert _{E_{2}}dy\leq C\int\limits_{\Omega }\left\Vert
\sum\limits_{j=1}^{m}r_{j}\left( y\right) u_{j}\right\Vert _{E_{1}}dy,
\end{equation*}%
where $\left\{ r_{j}\right\} $ is a sequence of independent symmetric $%
\left\{ -1,1\right\} $-valued random variables on $\Omega $. The smallest $C$
for which the above estimate holds is called a $R$-bound of the collection $%
G $ and denoted by $R\left( G\right) .$

Let $G_{h}$ be subset of $B\left( E_{1},E_{2}\right) $ depending on the
parameter $h\in Q.$ Here, $G_{h}$ is called uniform $R$-bounded in $h$ if
there is a constant $C$ independent on $h\in Q,$ such that 
\begin{equation*}
\sup\limits_{h\in Q}R\left( G_{h}\right) \leq C.
\end{equation*}

\textbf{Definition 1}$.$ A Banach space $E$\ is said to be a space
satisfying a multiplier condition if, for any $\Psi \in C^{\left( 1\right)
}\left( \mathbb{R};B\left( E\right) \right) $ the $R$-boundedness of the set 
\begin{equation*}
\left\{ \xi ^{j}\frac{d}{d\xi }\Psi \left( \xi \right) :\xi \in \mathbb{R}%
\backslash \left\{ 0\right\} ,\text{ }j=0,1\right\} 
\end{equation*}%
implies that $\Psi $ is a Fourier multiplier, i.e. $\Psi \in $ $%
M_{p}^{p}\left( E\right) $ for any $p\in \left( 1,\infty \right) .$

Note that $UMD$ spaces satisfies the multiplier condition (see e.g. $\left[ 
\text{5, 23}\right] $).

If 
\begin{equation*}
\text{ }\sup\limits_{h\in Q}R\left( \left\{ \left\vert \xi \right\vert
^{j}D^{j}\Psi _{h}\left( \xi \right) \text{: }\xi \in R\backslash \left\{
0\right\} ,j=0,1\right\} \right) \leq K
\end{equation*}%
then $\Psi _{h}$ is called a uniform collection of Fourier multipliers.

The $\varphi $-positive operator $A$ is said to be $R$-positive in a Banach
space $E$ if the set 
\begin{equation*}
\left\{ \xi \left( A+\xi \right) ^{-1}\text{: }\xi \in S_{\varphi }\right\} 
\text{, }0\leq \varphi <\pi
\end{equation*}%
is $R$-bounded.

Let $E_{0}$ and $E$ be two Banach spaces. $E_{0}$ is continuously and
densely embedded into $E$. Let $m$ be a positive integer.\ $W_{p}^{m}\left(
a,b;E_{0},E\right) $ denotes the collection of $E$-valued functions $u\in
L_{p}\left( a,b;E_{0}\right) $ that have the generalized derivatives $%
u^{\left( m\right) }\in L_{p}\left( a,b;E\right) $ with the norm 
\begin{equation*}
\ \left\Vert u\right\Vert _{W_{p}^{l}\left( a,b;E_{0},E\right) }=\left\Vert
u\right\Vert _{L_{p}\left( a,b;E_{0}\right) }+\left\Vert u^{\left( m\right)
}\right\Vert _{L_{p}\left( a,b;E\right) }<\infty .
\end{equation*}%
\ \ $\ \ \ \ $For $E_{0}=E$ it denotes by $W_{p}^{m}\left( \Omega ;E\right)
. $

Let $\varepsilon \in \left( 0,\right. \left. \varepsilon _{0}\right] $ be a
parameter for some positive bounded numbers $\varepsilon _{0}.$We define in $%
W_{p}^{m}\left( a,b;E_{0},E\right) $ the following parameterized norm 
\begin{equation*}
\left\Vert u\right\Vert _{W_{p,\varepsilon }^{m}\left( a,b;E_{0},E\right)
}=\left\Vert u\right\Vert _{L_{p}\left( a,b;E_{0}\right) }+\left\Vert
\varepsilon u^{\left( m\right) }\right\Vert _{L_{p}\left( a,b;E\right) }.
\end{equation*}

From $\left[ 20\right] $ we obtain:

\textbf{Theorem A}$_{1}.$ Assume the following conditions are satisfied:

(1) $E$\ is a Banach space satisfying the uniform multiplier condition for $%
p\in \left( 1,\infty \right) $;

(2) $0\leq \mu \leq 1-\frac{j}{m}$, $j=1,2,...,m-1;$\ 

(3) $A$ is an $R$-positive operator in $E$ with\ $0\leq \varphi <\pi $.

Then:

(a) the embedding%
\begin{equation*}
D^{j}W_{p}^{m}\left( a,b;E\left( A\right) ,E\right) \subset L_{p}\left(
a,b;E\left( A^{1-\frac{j}{m}-\mu }\right) \right)
\end{equation*}
is continuous and there exists a positive constant $C_{\mu }$ such that

\begin{equation*}
\varepsilon ^{\frac{j}{m}}\left\Vert u^{\left( j\right) }\right\Vert
_{L_{p}\left( \Omega ;E\left( A^{1-\frac{j}{m}-\mu }\right) \right) }\leq
\end{equation*}%
\begin{equation*}
C_{\mu }\left[ h^{\mu }\left\Vert u\right\Vert _{W_{p,\varepsilon
}^{m}\left( a,b;E\left( A\right) ,E\right) }+h^{-\left( 1-\mu \right)
}\left\Vert u\right\Vert _{L_{p}\left( a,b;E\right) }\right]
\end{equation*}%
for all $u\in W_{p}^{m}\left( a,b;E\left( A\right) ,E\right) ;$

(b) If $A^{-1}\in \sigma _{\infty }\left( E\right) $ and $0<\mu \leq 1-\frac{%
j}{m}$ then the embedding 
\begin{equation*}
D^{j}W_{p}^{m}\left( a,b;E\left( A\right) ,E\right) \subset L_{p}\left(
a,b;E\left( A^{1-\frac{j}{m}-\mu }\right) \right)
\end{equation*}%
is compact.

\textbf{Theorem A}$_{2}.$ Suppose all conditions of Theorem A$_{1}$
satisfied and $0<\mu <1-\frac{j}{m}.$ Then the embedding 
\begin{equation*}
D^{j}W_{p}^{m}\left( a,b;E\left( A\right) ,E\right) \subset L_{p}\left(
a,b;\left( E\left( A\right) ,E\right) _{\frac{j}{m},p}\right)
\end{equation*}%
is continuous and there exists a positive constant $C_{\mu }$ such that for
all $u\in W_{p}^{m}\left( a,b;E\left( A\right) ,E\right) $ the uniform
estimate holds

\begin{equation*}
\varepsilon ^{\frac{j}{m}}\left\Vert u^{\left( j\right) }\right\Vert
_{L_{p}\left( a,b;\left( E\left( A\right) ,E\right) _{\frac{j}{m}+\mu
,p}\right) }\leq
\end{equation*}

\begin{equation*}
C_{\mu }\left[ h^{\mu }\left\Vert u\right\Vert _{W_{p,\varepsilon
}^{m}\left( a,b;E\left( A\right) ,E\right) }+h^{-\left( 1-\mu \right)
}\left\Vert u\right\Vert _{L_{p}\left( a,b;E\right) }\right] .
\end{equation*}%
In a similar way as $\left[ \text{22},\text{ \S 1.7.7, Theorem 2}\right] $
and $\left[ \text{24, \S\ 10.1}\right] $ we obtain, respectively:

\textbf{Theorem A}$_{3}.$ Let $m$, $j$ be integer numbers, $0\leq j\leq m-1, 
$ $\theta _{j}=\frac{pj+1}{pm}$ and $x_{0}\in \left[ 0,b\right] .$

Then the transformation $u\rightarrow u^{\left( j\right) }\left(
x_{0}\right) $ is bounded linear from $W_{p}^{m}\left( 0,b;E_{0},E\right) $
onto $\left( E_{0},E\right) _{\theta _{j},p}$ and the inequality holds 
\begin{equation*}
\varepsilon ^{\theta _{j}}\left\Vert u^{\left( j\right) }\left( x_{0}\right)
\right\Vert _{\left( E_{0},E\right) _{\theta _{j},p}}\leq C\left( \left\Vert
\varepsilon u^{\left( m\right) }\right\Vert _{L_{p}\left( 0,b;E\right)
}+\left\Vert u\right\Vert _{L_{p\ }\left( 0,b;E_{0}\right) }\right) .
\end{equation*}

\textbf{Theorem A}$_{4}.$ Let $m$, $j$ be integer numbers, $0\leq j\leq m-1$%
, $\theta _{j}=\frac{pj+1}{pm}$ and $x_{0}\in \left[ 0,b\right] .$

Then the transformation $u\rightarrow u^{\left( j\right) }\left(
x_{0}\right) $ is bounded linear from $W_{p}^{m}\left( 0,b;E\right) $ into $%
E $ and the following inequality holds

\begin{equation*}
\varepsilon ^{\theta _{j}}\left\Vert u^{\left( j\right) }\left( x_{0}\right)
\right\Vert _{E}\leq C\left( h^{1-\theta _{j}}\left\Vert tu^{\left( m\right)
}\right\Vert _{L_{p}\left( 0,b;E\right) }+h^{-\theta _{j}}\left\Vert
u\right\Vert _{L_{p\ }\left( 0,b;E\right) }\right) .
\end{equation*}

From $\left[ \text{4, Theorem 2.1}\right] $ we obtain

\textbf{Theorem A}$_{5}.$ Let $E$ be a Banach space, $A$ be a $\varphi $%
-positive operator in $E$ with bound $M,$ $0\leq \varphi <\pi .$ Let $m$ be
a positive integer, $p\in \left( 1,\infty \right) $ and $\alpha \in \left( 
\frac{1}{2p}\text{, }\frac{1}{2p}+m\right) .$ Then, for $\lambda \in
S_{\varphi }$ the operator $-A_{\lambda }^{\frac{1}{2}}$ generates a
semigroup $e^{-xA_{\lambda }^{\frac{1}{2}}}$ which is holomorphic for $x>0.$
Moreover, there exists\ a positive constant $C$ (depending only on $%
M,\varphi ,m,\alpha $ and $p$) such that for every $u\in \left( E,E\left(
A^{m}\right) \right) _{\frac{\alpha }{m}-\frac{1}{2mp},p}$ and $\lambda \in
S_{\varphi },$%
\begin{equation}
\int_{0}^{\infty }\left\Vert A_{\lambda }^{\alpha }e^{-xA_{\lambda }^{\frac{1%
}{2}}}u\right\Vert ^{p}dx\leq M_{0}\left[ \left\Vert u\right\Vert _{\left(
E,E\left( A^{m}\right) \right) _{\frac{\alpha }{m}-\frac{1}{2mp}%
,p}}^{p}+\left\vert \lambda \right\vert ^{\alpha p-\frac{1}{2}}\left\Vert
u\right\Vert _{E}^{p}\right] .  \tag{1.3}
\end{equation}

Consider the nonlocal BVP for parameter dependent differential
operator-equation

\begin{equation*}
-\varepsilon u^{\left( 2\right) }\left( x,\varepsilon \right) +\left(
A+\lambda \right) u\left( x,\varepsilon \right) =0,
\end{equation*}%
\begin{equation*}
\sum\limits_{i=0}^{m_{k}}\varepsilon ^{\sigma _{i}}\left[ \alpha
_{ki}u^{\left( i\right) }\left( 0,\varepsilon \right) +\beta _{ki}u^{\left(
i\right) }\left( 1,\varepsilon \right) \right] =f_{k},k=1,2,
\end{equation*}%
where $f_{k}\in E,$ $\sigma _{i}=\frac{i}{2}+\frac{1}{2p},$ $p\in \left(
1,\infty \right) ,$ $m_{k}\in \left\{ 0,1\right\} ;$ $\alpha _{ki},$ $\beta
_{ki}$ are complex numbers; $\varepsilon $ is a positive and $\lambda $ is a
complex parameter; $A$ is a linear operator in $E.$ Let 
\begin{equation*}
E_{k}=\left( E\left( A\right) ,E\right) _{\theta _{k},p}\text{, }\theta _{k}=%
\frac{m_{k}}{2}+\frac{1}{2p}.
\end{equation*}

\textbf{Condition 1. }Let $\alpha _{k}=\alpha _{k,m_{k}},$ $\beta _{k}=\beta
_{k,m_{k}}.$ Suppose 
\begin{equation*}
d=\ \left( -1\right) ^{m_{1}}\alpha _{1}\beta _{2}-\left( -1\right)
^{m_{2}}\alpha _{2}\beta _{1}\neq 0,
\end{equation*}

and%
\begin{equation*}
\dsum\limits_{j-1}^{2}\dsum\limits_{i=0}^{2}\left\vert \alpha _{1i}\alpha
_{2j}\right\vert +\left\vert \beta _{1i}\beta _{2j}\right\vert <\left\vert
d\right\vert .
\end{equation*}

From $\left[ \text{17, Theorem 2 }\right] $ we obtain

\textbf{Theorem A}$_{6}$\textbf{. }Let the Condition 1 hold and $%
0<\varepsilon \leq \varepsilon _{0}$. Assume $E$\ is a Banach space
satisfying the uniform multiplier condition for $p\in \left( 1,\infty
\right) $ and\ $A$ is a $R$-positive operator in $E$ for $0\leq \varphi <\pi
.$Then problem $\left( 1.3\right) $ has a unique solution $u\in $\ $%
W_{p}^{2}\left( 0,1;E\left( A\right) ,E\right) $ for $f_{k}\in E_{k},$ $%
\theta _{k}=\frac{m_{k}}{2}+\frac{1}{2p},$ $p\in \left( 1,\infty \right) $, $%
\lambda \in S_{\varphi }$ with large enough $\left\vert \lambda \right\vert $
and the coercive uniform estimate holds 
\begin{equation*}
\sum\limits_{i=0}^{2}\varepsilon ^{\frac{i}{2}}\left\vert \lambda
\right\vert ^{1-\frac{i}{2}}\left\Vert u^{\left( i\right) }\right\Vert
_{L_{p}\left( 0,1;E\right) }+\left\Vert Au\right\Vert _{L_{p}\left(
0,1;E\right) }\leq M\sum\limits_{k=1}^{2}\left( \left\Vert f_{k}\right\Vert
_{E_{k}}+\left\vert \lambda \right\vert ^{1-\theta _{k}}\left\Vert
f_{k}\right\Vert _{E}\right) .
\end{equation*}

\begin{center}
\textbf{2. Abstract elliptic equation with parameters}
\end{center}

Consider the BVP for DOE with parameters

\begin{equation}
\left( L_{\varepsilon }+\lambda \right) u=\ -\varepsilon u^{\left( 2\right)
}\left( x,\varepsilon \right) +Au\left( x,\varepsilon \right) +Bu^{\left(
1\right) }\left( x,\varepsilon \right) +\lambda u\left( x,\varepsilon
\right) =f\left( x\right) ,\text{ }x\in \left( 0,T\right) ,  \tag{2.1}
\end{equation}%
\qquad \qquad \qquad

\begin{equation}
L_{1}u=\sum\limits_{i=0}^{m_{1}}\varepsilon ^{\frac{i}{2}}\alpha
_{i}u^{\left( i\right) }\left( 0,\varepsilon \right) =f_{1}\left(
\varepsilon \right) ,\text{ }L_{2}u=\sum\limits_{i=0}^{m_{2}}\varepsilon ^{%
\frac{i}{2}}\beta _{i}u^{\left( i\right) }\left( T,\varepsilon \right)
=f_{2}\left( \varepsilon \right) ,  \tag{2.2}
\end{equation}%
where $m_{k}\in \left\{ 0,1\right\} ,$ $\alpha _{i},$ $\beta _{i}$ are
complex numbers; $\varepsilon $ is a positive and $\lambda $ is a complex
parameter; $A$ and $B$ are linear operators in $E$ and $u\left( x\right)
=u\left( x,\varepsilon \right) $ is a solution of $\left( 2.1\right) -\left(
2.2\right) .$

\bigskip First all of, consider the problem $\left( 2.1\right) -\left(
2.2\right) $ with $f_{k}=0$, i.e. consider the homogenous problem%
\begin{equation}
\ -\varepsilon u^{\left( 2\right) }\left( x,\varepsilon \right) +Bu^{\left(
1\right) }\left( x,\varepsilon \right) +\left( A+\lambda \right) u\left(
x,\varepsilon \right) =0,\text{ }x\in \left( 0,T\right) ,  \tag{2.3}
\end{equation}%
\begin{equation}
\sum\limits_{i=0}^{m_{k}}\varepsilon ^{\frac{i}{2}}\alpha _{i}u^{\left(
i\right) }\left( 0,\varepsilon \right) =f_{1}\left( \varepsilon \right) ,%
\text{ }\sum\limits_{i=0}^{m_{k}}\varepsilon ^{\frac{i}{2}}\beta
_{i}u^{\left( i\right) }\left( T,\varepsilon \right) =f_{2}\left(
\varepsilon \right) ,  \tag{2.4}
\end{equation}%
where%
\begin{equation*}
f_{k}=f_{k}\left( \varepsilon \right) \in E_{k}=\left( E\left( A\right)
,E\right) _{\theta _{k},p}\text{ for all }\varepsilon >0,\text{ }
\end{equation*}

\begin{equation*}
\theta _{k}=\frac{m_{k}}{2}+\frac{1}{2p}\text{, }m_{k}\in \left\{
0,1\right\} ,\text{ }k=1,\text{ }2\text{, }p\in \left( 1,\infty \right) .
\end{equation*}

Let 
\begin{equation*}
d=\alpha _{0}\beta _{1}-\beta _{0}\alpha _{1}\text{, }X=L_{p}\left(
0,T;E\right) \text{, }Y=\ W_{p}^{2}\left( 0,T;E\left( A\right) ,E\right) .
\end{equation*}

\textbf{Condition 2.1.} Assume the following conditions are satisfied:

(1) Assume $E$\ is a Banach space satisfying the uniform multiplier
condition for $p\in \left( 1,\infty \right) $;

(2) $A$ is a $R$-positive operator in $E$ for $0\leq \varphi <\pi $ and $%
d\neq 0$;

(3) $B$ is a bounded operator, $\left( A+B\right) ^{-\frac{1}{2}}\in B\left(
E\right) $ and 
\begin{equation*}
\left\Vert B\right\Vert _{B\left( E\right) }<\sup\limits_{t\in \left[
0,\infty \right] }\left\Vert A\left( A+t\right) ^{-1}\right\Vert _{B\left(
E\right) }.
\end{equation*}

\textbf{Theorem 2.1. }Assume the Condition 2.1 hold. Then problem $\left(
2.3\right) -\left( 2.4\right) $ has a unique solution $u\in $\ $Y$ for $%
f_{k}\in E_{k},$ $\lambda \in S_{\varphi }$ with large enough $\left\vert
\lambda \right\vert .$ Moreover, the coercive estimate holds

\begin{equation}
\sum\limits_{i=0}^{2}\varepsilon ^{\frac{i}{2}-\frac{1}{p}}\left\vert
\lambda \right\vert ^{1-\frac{i}{2}}\left\Vert u^{\left( i\right) }\left(
.,\varepsilon \right) \right\Vert _{X}+\left\Vert Au\right\Vert _{X}\leq
M\sum\limits_{k=1}^{2}\left( \left\Vert f_{k}\right\Vert _{E_{k}}+\left\vert
\lambda \right\vert ^{1-\theta _{k}}\left\Vert f_{k}\right\Vert _{E}\right) 
\tag{2.5}
\end{equation}%
uniformly with respect to $\varepsilon $ and $\lambda .$\qquad \qquad \qquad

\textbf{Proof: } By definition of positive operator, $4\varepsilon A$ is $%
\varphi $-positive uniformly in $\varepsilon \in \left( 0\right. ,\left. 1%
\right] .$ Then for $\left\vert \arg \lambda \right\vert \leq \varphi $, $%
\left\vert \arg \mu \right\vert \leq \varphi _{1}$ and $\varphi +\varphi
_{1}<\pi $ we have the estimate \bigskip

\begin{equation*}
\left\Vert \left( 4\varepsilon A_{\lambda }+\mu \right) ^{-1}\right\Vert
\leq \frac{M_{0}}{\left\vert \mu \right\vert }
\end{equation*}%
where $A_{\lambda }=A+\lambda $ and $M_{0}$ depend only on $\varphi $. By
perturbation theory of positive operators and semigroups (see e.g. $\left[ 
\text{14, \S\ 1.3}\right] $ and $\left[ \text{7, \S\ 3}\right] $) there
exists the analytic semigroups 
\begin{equation*}
U_{\lambda }\left( x,\varepsilon \right) =\exp -\left\{ \varepsilon
^{-1}xA_{\lambda }^{\frac{1}{2}}\right\} .
\end{equation*}
Moreover, by virtue of Condition 2.1 and in view of the same perturbation
theory, the following semigroups 
\begin{equation*}
U_{1,\lambda }\left( x,\varepsilon \right) =\exp -\left\{ xQ_{1,\lambda
}\left( \varepsilon \right) \right\} ,\text{ }U_{2,\lambda }\left(
x,\varepsilon \right) =\exp -\left\{ xQ_{2,\lambda }\left( \varepsilon
\right) \right\}
\end{equation*}%
are holomorphic for $x>0$ and strongly continuous for $x\geq 0,$ where%
\begin{equation}
Q_{1,\lambda }\left( \varepsilon \right) =\frac{1}{2\varepsilon }\left[
B+\left( B^{2}+4\varepsilon A_{\lambda }\right) ^{\frac{1}{2}}\right] ,\text{
}Q_{2,\lambda }\left( \varepsilon \right) =\frac{1}{2\varepsilon }\left[
B-\left( B^{2}+4\varepsilon A_{\lambda }\right) ^{\frac{1}{2}}\right] . 
\tag{2.6 }
\end{equation}

\bigskip Let firstly, show that the function $u\left( x,\varepsilon \right)
=U_{1,\lambda }\left( x,\varepsilon \right) g_{1}+U_{2,\lambda }\left(
x,\varepsilon \right) g_{2}$ is a solution of the equation $\left(
2.3\right) $ belonging $Y$ for 
\begin{equation*}
g_{1},\text{ }g_{2}\in \left( E\left( A\right) ,E\right) _{\frac{1}{2p},p%
\text{ }}.
\end{equation*}%
Indeed, by properties of continuous semigroups it is clear to see that
operator functions $U_{1\varepsilon }\left( x\right) $ and $U_{2\varepsilon
}\left( x\right) $ are solution of $\left( 2.3\right) .$ From $\left(
2.6\right) $ we get 
\begin{equation*}
\frac{d^{2}u}{dx^{2}}=Q_{1,\lambda }^{2}\left( \varepsilon \right)
U_{1,\lambda }\left( x,\varepsilon \right) g_{1}+Q_{2,\lambda }^{2}\left(
\varepsilon \right) U_{2,\lambda }\left( x,\varepsilon \right) g_{2},\text{ }
\end{equation*}%
\begin{equation*}
Au\left( x,\varepsilon \right) =A\left[ U_{1,\lambda }\left( x,\varepsilon
\right) g_{1}+U_{2,\lambda }\left( x,\varepsilon \right) g_{2}\right] .
\end{equation*}

Then 
\begin{equation*}
\left\Vert u\right\Vert _{Y}=\left\Vert Au\right\Vert _{X}+\left\Vert
u^{\left( 2\right) }\right\Vert _{X}\leq
\end{equation*}

\begin{equation}
\left( \dint\limits_{0}^{T}\left\Vert AU_{1,\lambda }\left( x,\varepsilon
\right) g_{1}\right\Vert _{E}^{p}dx\right) ^{\frac{1}{p}}+\left(
\dint\limits_{0}^{T}\left\Vert AU_{2,\lambda }\left( x,\varepsilon \right)
g_{2}\right\Vert _{E}^{p}dx\right) ^{\frac{1}{p}}+  \tag{2.7}
\end{equation}%
\begin{equation*}
\left( \dint\limits_{0}^{T}\left\Vert Q_{1,\lambda }^{2}\left( \varepsilon
,\lambda \right) U_{1,\lambda }\left( x,\varepsilon \right) g_{1}\right\Vert
_{E}^{p}dx\right) ^{\frac{1}{p}}+\left( \dint\limits_{0}^{T}\left\Vert
Q_{2,\lambda }^{2}\left( \varepsilon \right) U_{2,\lambda }\left(
x,\varepsilon \right) g_{2}\right\Vert _{E}^{p}dx\right) ^{\frac{1}{p}}.
\end{equation*}

By properties of positive operators and by Theorem A$_{5}$ we have%
\begin{equation}
\left( \dint\limits_{0}^{T}\left\Vert AU_{1,\lambda }\left( x,\varepsilon
\right) g_{1}\right\Vert _{E}^{p}dx\right) ^{\frac{1}{p}}+\left(
\dint\limits_{0}^{T}\left\Vert AU_{2,\lambda }\left( x,\varepsilon \right)
g_{2}\right\Vert _{E}^{p}dx\right) ^{\frac{1}{p}}\leq  \tag{2.8}
\end{equation}%
\begin{equation*}
\left( 1+\left\Vert AA_{\lambda }^{-1}\right\Vert _{B\left( E\right)
}\right) \left[ \left( \dint\limits_{0}^{T}\left\Vert A_{\lambda
}U_{1,\lambda }\left( x,\varepsilon \right) g_{1}\right\Vert
_{E}^{p}dx\right) ^{\frac{1}{p}}\right. +
\end{equation*}%
\begin{equation*}
\left. \left( \dint\limits_{0}^{T}\left\Vert A_{\lambda }U_{2,\lambda
}\left( x,\varepsilon \right) g_{2}\right\Vert _{E}^{p}dx\right) ^{\frac{1}{p%
}}\right] \leq C_{0}\left\Vert U_{1,\lambda }\left( x,\varepsilon \right)
V_{\lambda }^{-1}\left( x\right) \right\Vert _{B\left( E\right) }
\end{equation*}%
\begin{equation*}
\left[ \left( \dint\limits_{0}^{T}\left\Vert A_{\lambda }V_{\lambda }\left(
x\right) g_{1}\right\Vert _{E}^{p}dx\right) ^{\frac{1}{p}}+\left(
\dint\limits_{0}^{T}\left\Vert A_{\lambda }V_{\lambda }\left( x\right)
g_{2}\right\Vert _{E}^{p}dx\right) ^{\frac{1}{p}}\right] \leq
\end{equation*}%
\begin{equation*}
C_{0}N_{0}M_{0}\dsum\limits_{k=1}^{2}\left( \left\Vert g_{k}\right\Vert
_{\left( E\left( A\right) ,E\right) _{\frac{1}{2p},p}}+\left\vert \lambda
\right\vert ^{1-\frac{1}{2p}}\left\Vert g_{k}\right\Vert _{E}\right) ,
\end{equation*}%
where $M_{0}$ is a constant in $\left( 1.3\right) $ and 
\begin{equation*}
C_{0}=\left( 1+\left\Vert AA_{\lambda }^{-1}\right\Vert _{B\left( E\right)
}\right) \text{, }N_{0}=\left\Vert U_{1,\lambda }\left( x,\varepsilon
\right) V_{\lambda }^{-1}\left( x\right) \right\Vert _{B\left( E\right) }%
\text{ for }\lambda \in S\left( \varphi \right) .
\end{equation*}

In a similar way, we get the uniform estimate 
\begin{equation*}
\left( \dint\limits_{0}^{T}\left\Vert Q_{1,\lambda }^{2}\left( \varepsilon
,\lambda \right) U_{1,\lambda }\left( x,\varepsilon \right) g_{1}\right\Vert
_{E}^{p}dx\right) ^{\frac{1}{p}}+\left( \dint\limits_{0}^{T}\left\Vert
Q_{2,\lambda }^{2}\left( \varepsilon \right) U_{2,\lambda }\left(
x,\varepsilon \right) g_{2}\right\Vert _{E}^{p}dx\right) ^{\frac{1}{p}}\leq
\end{equation*}

\begin{equation}
M_{1}\dsum\limits_{k=1}^{2}\left( \left\Vert g_{k}\right\Vert _{\left(
E\left( A\right) ,E\right) _{\frac{1}{2p},p}}+\left\vert \lambda \right\vert
^{1-\frac{1}{2p}}\left\Vert g_{k}\right\Vert _{E}\right) .  \tag{2.9}
\end{equation}

From $\left( 2.7\right) ,$ $\left( 2.8\right) $ and $\left( 2.9\right) $ we
obtain that 
\begin{equation*}
u\left( \varepsilon ,.\right) \in W_{p}^{2}\left( 0,T;E\left( A\right)
,E\right) \text{ for }g_{1},g_{2}\in \left( E\left( A\right) ,E\right) _{%
\frac{1}{2p},p\text{ }}.
\end{equation*}%
Without loss of generality assume $m_{1}=m_{2}=1.$ A function 
\begin{equation*}
u\left( x,\varepsilon \right) =U_{1,\lambda }\left( x,\varepsilon \right)
g_{1}+U_{2,\lambda }\left( x,\varepsilon \right) g_{2}
\end{equation*}%
satisfies the boundary conditions $\left( 2.4\right) $ if 
\begin{equation}
\left( \varepsilon \alpha _{1}Q_{1,\lambda }\left( \varepsilon \right)
+\alpha _{0}\right) g_{1}+\left( \varepsilon \alpha _{1}Q_{2,\lambda }\left(
\varepsilon \right) +\alpha _{0}\right) g_{2}=f_{1},  \tag{2.10}
\end{equation}

\begin{equation*}
\left( \varepsilon \beta _{1}Q_{1,\lambda }\left( \varepsilon \right) +\beta
_{0}\right) g_{1}+\left( \varepsilon \beta _{1}Q_{2,\lambda }\left(
\varepsilon \right) +\beta _{0}\right) g_{2}=f_{2}.
\end{equation*}

The main operator-determinant of the algebraic equation $\left( 2.10\right) $
(with respect to $g_{1}$ and $g_{2}$) can be expressed as 
\begin{equation*}
D_{\lambda }\left( \varepsilon \right) =\varepsilon ^{2}\alpha _{1}\beta
_{1}Q_{1,\lambda }\left( \varepsilon \right) Q_{2,\lambda }\left(
\varepsilon \right) +\varepsilon \alpha _{1}\beta _{0}Q_{1,\lambda }\left(
\varepsilon \right) +\alpha _{0}\varepsilon \beta _{1}Q_{2,\lambda }\left(
\varepsilon \right) +\alpha _{0}\beta _{0}-
\end{equation*}%
\begin{equation*}
\varepsilon ^{2}\alpha _{1}\beta _{1}Q_{1,\lambda }\left( \varepsilon
\right) Q_{2,\lambda }\left( \varepsilon \right) +\varepsilon \alpha
_{0}\beta _{1}Q_{1,\lambda }\left( \varepsilon \right) +\beta
_{0}\varepsilon \alpha _{1}Q_{2,\lambda }\left( \varepsilon \right) +\alpha
_{0}\beta _{0}=
\end{equation*}%
\begin{equation*}
\varepsilon \left( \alpha _{1}\beta _{0}-\alpha _{0}\beta _{1}\right)
Q_{1,\lambda }\left( \varepsilon \right) +\varepsilon \left( \alpha
_{0}\beta _{1}-\beta _{0}\alpha _{1}\right) Q_{2,\lambda }\left( \varepsilon
\right) =\varepsilon d\left[ Q_{2,\lambda }\left( \varepsilon \right)
-Q_{1,\lambda }\left( \varepsilon \right) \right] .
\end{equation*}

Since $d\neq 0,$ 
\begin{equation*}
\left[ Q_{2,\lambda }\left( \varepsilon \right) -Q_{1,\lambda }\left(
\varepsilon \right) \right] =\frac{1}{\varepsilon }Q_{\lambda }\left(
\varepsilon \right) =\frac{1}{\varepsilon }\left( B^{2}+4\varepsilon
A_{\lambda }\right) ^{\frac{1}{2}},\text{ }
\end{equation*}%
where 
\begin{equation*}
Q_{\lambda }\left( \varepsilon \right) =\left( B^{2}+4\varepsilon A_{\lambda
}\right) ^{\frac{1}{2}}.
\end{equation*}%
It is clear to see that $Q_{\lambda }\left( \varepsilon \right) $ has a
bounded inverse $Q_{\lambda }^{-1}\left( \varepsilon \right) $. Hence, $%
D_{\lambda }\left( \varepsilon \right) $ has a bounded inverse 
\begin{equation}
D_{\lambda }^{-1}\left( \varepsilon \right) =-d^{-1}Q_{\lambda }^{-1}\left(
\varepsilon \right)  \tag{2.11}
\end{equation}%
for $\varepsilon >0$ and $\lambda \in S\left( \varphi \right) .$ So, the
system $\left( 2.10\right) $ has a unique solution%
\begin{equation}
g_{1}=D_{1,\lambda }\left( \varepsilon \right) D_{\lambda }^{-1}\left(
\varepsilon \right) ,\text{ }g_{2}=D_{2,\lambda }\left( \varepsilon \right)
D_{\lambda }^{-1}\left( \varepsilon \right) ,  \tag{2.12}
\end{equation}%
where 
\begin{equation*}
D_{1,\lambda }\left( \varepsilon \right) =\left\vert 
\begin{array}{cc}
f_{1} & \varepsilon \alpha _{1}Q_{2,\lambda }\left( \varepsilon \right)
+\alpha _{0} \\ 
f_{2} & \varepsilon \beta _{1}Q_{2,\lambda }\left( \varepsilon \right)
U_{2,\lambda }\left( 1,\varepsilon \right) +\beta _{0}U_{2,\lambda }\left(
1,\varepsilon \right)%
\end{array}%
\right\vert =
\end{equation*}%
\begin{equation*}
\left[ \varepsilon \beta _{1}Q_{2,\lambda }\left( \varepsilon \right)
U_{2,\lambda }\left( 1,\varepsilon \right) +\beta _{0}U_{2,\lambda }\left(
1,\varepsilon \right) \right] f_{1}-\left[ \varepsilon \alpha
_{1}Q_{2,\lambda }\left( \varepsilon \right) +\alpha _{0}\right] f_{2},
\end{equation*}

\begin{equation*}
D_{2,\lambda }\left( \varepsilon \right) =\left\vert \left[ 
\begin{array}{cc}
\varepsilon \alpha _{1}Q_{1,\lambda }\left( \varepsilon \right) +\alpha _{0}
& f_{1} \\ 
\varepsilon \beta _{1}Q_{1,\lambda }\left( \varepsilon \right) U_{1,\lambda
}\left( \varepsilon ,1\right) +\beta _{0}U_{1,\lambda }\left( \varepsilon
,1\right) & f_{2}%
\end{array}%
\right] \right\vert =
\end{equation*}%
\begin{equation*}
\left[ \varepsilon \alpha _{1}Q_{1,\lambda }\left( \varepsilon \right)
+\alpha _{0}\right] f_{2}-\left[ \varepsilon \beta _{1}Q_{1,\lambda }\left(
\varepsilon \right) U_{1,\lambda }\left( \varepsilon ,1\right) +\beta
_{0}U_{1,\lambda }\left( 1,\varepsilon \right) \right] f_{1}.
\end{equation*}

From $\left( 2.7\right) $ and $\left( 2.11\right) $ we get the following
representation of solution $\left( 2.3\right) -\left( 2.4\right) :$

\begin{equation}
u\left( x,\varepsilon \right) =D_{\lambda }^{-1}\left( \varepsilon \right) 
\left[ U_{1,\lambda }\left( x,\varepsilon \right) D_{1,\lambda }\left(
\varepsilon \right) +U_{2,\lambda }\left( x,\varepsilon \right) D_{2,\lambda
}\left( \varepsilon \right) \right] =  \tag{2.13}
\end{equation}

\begin{equation*}
D_{\lambda }^{-1}\left( \varepsilon \right) \left\{ U_{1,\lambda }\left(
x,\varepsilon \right) U_{2,\lambda }\left( 1,\varepsilon \right) \left[
\varepsilon \beta _{1}Q_{2,\lambda }\left( \varepsilon \right) +\beta _{0}%
\right] \right. -
\end{equation*}%
\begin{equation*}
U_{2,\lambda }\left( x,\varepsilon \right) U_{1,\lambda }\left(
1,\varepsilon \right) \left[ \varepsilon \beta _{1}Q_{1,\lambda }\left(
\varepsilon \right) +\beta _{0}\right] f_{1}+
\end{equation*}%
\begin{equation*}
D_{\lambda }^{-1}\left( \varepsilon \right) \left\{ U_{2,\lambda }\left(
x,\varepsilon \right) \left[ \varepsilon \alpha _{1}Q_{1,\lambda }\left(
\varepsilon \right) +\alpha _{0}\right] -\right. \left. U_{1,\lambda }\left(
x,\varepsilon \right) \left[ \varepsilon \alpha _{1}Q_{2,\lambda }\left(
\varepsilon \right) +\alpha _{0}\right] f_{2}\right\} .
\end{equation*}

Due to uniform boundedness of $D_{\lambda }^{-1}\left( \varepsilon \right) $
from $\left( 2.7\right) $ we obtain 
\begin{equation}
\sum\limits_{i=0}^{2}\varepsilon ^{\frac{i}{2}}\left\vert \lambda
\right\vert ^{1-\frac{i}{2}}\left\Vert u^{\left( i\right) }\right\Vert
_{X}+\left\Vert Au\right\Vert _{X}\leq C\sum\limits_{i=0}^{2}\varepsilon ^{%
\frac{i}{2}}\left\vert \lambda \right\vert ^{1-\frac{i}{2}}  \tag{2.14}
\end{equation}%
\begin{equation*}
\left\{ \dsum\limits_{k=1}^{2}\right. \left\Vert \varepsilon U_{3-k,\lambda
}\left( 1,\varepsilon \right) Q_{3-k,\lambda }\left( \varepsilon \right)
Q_{k,\lambda }^{i}\left( \varepsilon \right) U_{k,\lambda }\left(
x,\varepsilon \right) f_{1}\right\Vert _{X}+
\end{equation*}%
\begin{equation*}
\dsum\limits_{k=1}^{2}\left\Vert U_{3-k,\lambda }\left( 1,\varepsilon
\right) Q_{k,\lambda }^{i}\left( \varepsilon \right) U_{k,\lambda }\left(
x,\varepsilon \right) f_{1}\right\Vert _{X}+\dsum\limits_{k=1}^{2}\left\Vert
\varepsilon Q_{3-k,\lambda }\left( \varepsilon \right) Q_{k,\lambda
}^{i}\left( \varepsilon \right) U_{k,\lambda }\left( x,\varepsilon \right)
f_{2}\right\Vert _{X}+
\end{equation*}

\begin{equation*}
\dsum\limits_{k=1}^{2}\left\Vert Q_{k,\lambda }^{i}\left( \varepsilon
\right) U_{k,\lambda }\left( x,\varepsilon \right) f_{2}\right\Vert
_{X}+\dsum\limits_{k=1}^{2}\left\Vert \varepsilon U_{3-k,\lambda }\left(
1,\varepsilon \right) Q_{3-k,\lambda }\left( \varepsilon \right)
AU_{k,\lambda }\left( x,\varepsilon \right) f_{1}\right\Vert _{X}+
\end{equation*}%
\begin{equation*}
\dsum\limits_{k=1}^{2}\left\Vert U_{3-k,\lambda }\left( 1,\varepsilon
\right) AU_{k,\lambda }\left( x,\varepsilon \right) f_{1}\right\Vert
_{X}+\dsum\limits_{k=1}^{2}\left\Vert \varepsilon Q_{3-k,\lambda }\left(
\varepsilon \right) AU_{k,\lambda }\left( x,\varepsilon \right)
f_{2}\right\Vert _{X}+
\end{equation*}%
\begin{equation*}
\left. \dsum\limits_{k=1}^{2}\left\Vert AU_{k,\lambda }\left( x,\varepsilon
\right) f_{2}\right\Vert _{X}\right\} .
\end{equation*}

By $\left[ \text{4, Lemma 2.6}\right] $, we have

\begin{equation}
\left\Vert A^{\alpha }A_{\lambda }^{\beta }\right\Vert _{B\left( E\right)
}\leq C\left( 1+\left\vert \lambda \right\vert \right) ^{\alpha -\beta },%
\text{ }0\leq \alpha \leq \beta ,  \tag{2.15}
\end{equation}%
\begin{equation*}
\left\Vert A_{\lambda }^{\alpha }U_{k,\lambda }\left( x,\varepsilon \right)
\right\Vert _{B\left( E\right) }\leq Ce^{-\omega \varepsilon
^{-1}x\left\vert \lambda \right\vert ^{\frac{1}{2}}},\text{ for }\alpha \in 
\mathbb{R}\text{, }x\geq x_{0}>0\text{, }\lambda \in S\left( \varphi \right)
.
\end{equation*}

By properties of positive operators, from $\left( 2.6\right) $ and $\left(
2.15\right) $ for $u\in D\left( A^{\frac{1}{2}}\right) $\ we get 
\begin{equation*}
\text{ }Q_{\lambda }\left( \varepsilon \right) =Q_{\lambda }\left(
\varepsilon \right) A_{\lambda }^{-\frac{1}{2}}A_{\lambda }^{\frac{1}{2}},%
\text{ }\left\Vert \text{ }Q_{\lambda }\left( \varepsilon \right)
u\right\Vert _{E}\leq
\end{equation*}%
\begin{equation}
\left\Vert \text{ }Q_{\lambda }\left( \varepsilon \right) A_{\lambda }^{-%
\frac{1}{2}}\right\Vert _{B\left( E\right) }\left\Vert A_{\lambda }^{\frac{1%
}{2}}u\right\Vert _{E}\leq C\left\Vert A_{\lambda }^{\frac{1}{2}%
}u\right\Vert _{E}.  \tag{2.16}
\end{equation}

Moreover, by virtue of analytic semigroups theory, for all $u\in E$\ we have 
\begin{equation*}
\left\Vert U_{k,\lambda }\left( x,\varepsilon \right) u\right\Vert _{E}\leq
C\left\Vert U_{\lambda }\left( x,\varepsilon \right) u\right\Vert _{E}\text{%
, }k=1\text{, }2.
\end{equation*}

\bigskip By chance of variable, by estimates $\left( 2.14\right) -\left(
2.16\right) $ and by virtue of Theorem 1.5 we obtain

\begin{equation*}
\sum\limits_{i=0}^{2}\varepsilon ^{\frac{i}{2}-\frac{1}{2p}}\left\vert
\lambda \right\vert ^{1-\frac{i}{2}}\left\Vert u^{\left( i\right)
}\right\Vert _{X}+\left\Vert Au\right\Vert _{X}\leq M_{1}\left\Vert
A_{\lambda }^{-\left( 1-\frac{m_{k}}{2}\right) }\right\Vert _{B\left(
E\right) }\sum\limits_{i=0}^{2}\varepsilon ^{\frac{i}{2}-\frac{1}{p}%
}\left\vert \lambda \right\vert ^{1-\frac{i}{2}}
\end{equation*}%
\begin{equation*}
\sum\limits_{k=1}^{2}\left\Vert A_{\lambda }^{\left( 1-\frac{m_{k}}{2}%
\right) }U_{\lambda }\left( x,\varepsilon \right) f_{k}\right\Vert _{X}\leq
M\sum\limits_{k=1}^{2}\left[ \left\Vert f_{k}\right\Vert _{E_{k}}+\left\vert
\lambda \right\vert ^{1-\theta _{k}}\left\Vert f_{k}\right\Vert \right] .
\end{equation*}

\bigskip \textbf{Remark 2.1. }It is clear to see that the solution of the
problem $\left( 2.3\right) -\left( 2.4\right) $ depends on $\varepsilon $,
i.e. $u=u\left( x,\varepsilon \right) .$ Hence, it is interesting to
investigate behavior of solution when $\varepsilon \rightarrow 0$ and to
have the smoothness properties of the solution with respect to parameter $%
\varepsilon .$ From Theorem 3.1 we obtain the following result

\textbf{Corollary 2.1. }Assume all conditions of Theorem 2.1 are satisfied.
Then the solution $u$ of the problem $\left( 2.3\right) -\left( 2.4\right) $
satisfies the following:

(a) $\varepsilon ^{\frac{1}{p}}u\left( x,\varepsilon \right) =O\left(
\sum\limits_{k=1}^{2}A_{\lambda }^{-\frac{m_{k}}{2}}f_{k}\right) $ when $%
\varepsilon \rightarrow 0;$

(b) 
\begin{equation}
\varepsilon ^{\frac{3}{2}-\frac{1}{p}}\left\vert \lambda \right\vert ^{\frac{%
1}{2}}\left\Vert \frac{du}{d\varepsilon }\right\Vert _{X}+\varepsilon ^{3-%
\frac{1}{p}}\left\Vert \frac{d^{2}u}{d\varepsilon ^{2}}\right\Vert _{X}\leq
C\sum\limits_{k=1}^{2}\left( \left\Vert f_{k}\right\Vert _{E_{k}}+\left\vert
\lambda \right\vert ^{1-\theta _{k}}\left\Vert f_{k}\right\Vert _{E}\right) .
\tag{2.17}
\end{equation}

\textbf{Proof. }The part (a) is obtained from the representation of solution 
$\left( 2.13\right) .$ By differentiating both parts of $\left( 2.13\right) $
with respect to $\varepsilon $ and by using Theorem 1.5, the part (b) is
obtained.

\textbf{Theorem 2.2. }Assume the Condition 2.1 hold.\ Then the operator $%
u\rightarrow \left\{ \left( L_{\varepsilon }+\lambda \right)
u,L_{1}u,L_{2}u\right\} $ is an isomorphism from\ $Y$ onto $X\times
E_{1}\times E_{2}$ for $\left\vert \arg \lambda \right\vert \leq \varphi $, $%
0\leq \varphi <\pi $ with large enough $\left\vert \lambda \right\vert $.
Moreover, the uniform coercive estimate holds:

\begin{center}
\begin{equation}
\sum\limits_{j=0}^{2}\varepsilon ^{\frac{j}{2}}\left\vert \lambda
\right\vert ^{1-\frac{j}{2}}\left\Vert u^{\left( j\right) }\right\Vert
_{X}+\left\Vert Au\right\Vert _{X}\leq C\left[ \left\Vert f\right\Vert
_{X}+\sum\limits_{k=1}^{2}\left( \left\Vert f_{k}\right\Vert
_{E_{k}}+\left\vert \lambda \right\vert ^{1-\theta _{k}}\left\Vert
f_{k}\right\Vert _{E}\right) \right] .  \tag{2.18}
\end{equation}
\end{center}

\textbf{Proof. }We have proved the uniqueness of solution of $\left(
2.1\right) -\left( 2.2\right) $ in Theorem 2.1. Let us define 
\begin{equation*}
\bar{f}\left( x\right) =\left\{ 
\begin{array}{c}
f\left( x\right) \text{ if }x\in \left[ 0,T\right] \\ 
0\text{ \ \ \ \ \ if \ }x\notin \left[ 0,T\right]%
\end{array}%
\right\} .
\end{equation*}%
We now show that problem $\left( 2.1\right) -\left( 2.2\right) $ has a
solution $u\in Y$\ for all $f\in X$, $f_{k}\in E_{k}$ and $u=u_{1}+u_{2}$,
where $u_{1}$ is the restriction on $\left[ 0,1\right] $ of the solution of
the equation

\begin{equation}
\left( L_{\varepsilon }+\lambda \right) u=\bar{f}\left( x\right) ,\text{ }%
x\in \mathbb{R}=\left( -\infty ,\infty \right)  \tag{2.19}
\end{equation}%
and $u_{2}$ is a solution of the problem

\begin{equation}
\left( L_{\varepsilon }+\lambda \right) u=0\text{, }L_{k}u=f_{k}-L_{k}u_{1}.
\tag{2.20}
\end{equation}

By applying the Fourier transform we get that, the solution $\left(
2.19\right) $ can be given by 
\begin{equation*}
u\left( x\right) =F^{-1}\Phi \left( \lambda ,\varepsilon ,\xi \right) F\bar{f%
}=\frac{1}{2\pi }\int\limits_{\infty }^{\infty }e^{i\xi x}\Phi \left(
\lambda ,\varepsilon ,\xi \right) \left( F\bar{f}\right) \left( \xi \right)
d\xi ,
\end{equation*}%
where 
\begin{equation*}
\Phi \left( \lambda ,\varepsilon ,\xi \right) =\left( A-i\xi B+\varepsilon
\xi ^{2}+\lambda \right) ^{-1},
\end{equation*}%
here $i$ is the complex unity. It follows from the above expression that

\begin{center}
\begin{equation}
\sum\limits_{j=0}^{2}\varepsilon ^{\frac{j}{2}}\left\vert \lambda
\right\vert ^{1-\frac{j}{2}}\left\Vert u^{\left( j\right) }\right\Vert
_{L_{p}\left( R;E\right) }+\left\Vert Au\right\Vert _{L_{p}\left( R;E\right)
}=  \tag{2.21}
\end{equation}%
\begin{equation*}
\sum\limits_{j=0}^{2}\varepsilon ^{\frac{j}{2}}\left\vert \lambda
\right\vert ^{1-\frac{j}{2}}\left\Vert F^{-1}\xi ^{j}\Phi \left( \lambda
,\varepsilon ,\xi \right) F\bar{f}\right\Vert _{L_{p}\left( R;E\right)
}+\left\Vert F^{-1}A\Phi \left( \lambda ,\varepsilon ,\xi \right) F\bar{f}%
\right\Vert _{L_{p}\left( R;E\right) }.
\end{equation*}
\end{center}

\begin{quote}
Let us show that operator-functions%
\begin{equation*}
\Psi \left( \lambda ,\varepsilon ,\xi \right) =A\Phi \left( \lambda
,\varepsilon ,\xi \right) ,\sigma \left( \lambda ,\varepsilon ,\xi \right)
=\sum\limits_{j=0}^{2}\varepsilon ^{\frac{j}{2}}\left\vert \lambda
\right\vert ^{1-\frac{j}{2}}\xi ^{j}\Phi \left( \lambda ,\varepsilon ,\xi
\right)
\end{equation*}
are Fourier multipliers in $L_{p}\left( R;E\right) $. Actually, due to
positivity of $A$ and by assumption (2)\ we have
\end{quote}

\begin{equation}
\left\Vert \Phi \left( \lambda ,\varepsilon ,\xi \right) \right\Vert
_{B\left( E\right) }\leq M\left( 1+\left\vert \varepsilon \xi ^{2}+\lambda
\right\vert \right) ^{-1}\leq C_{1},  \tag{2.22}
\end{equation}%
\begin{equation*}
\text{ }\left\Vert \Psi \left( \lambda ,\varepsilon ,\xi \right) \right\Vert
_{B\left( E\right) }=\left\Vert A\Phi \left( \lambda ,\varepsilon ,\xi
\right) \right\Vert \leq C_{2}.
\end{equation*}%
It is clear to observe that 
\begin{equation*}
\xi \frac{d}{d\xi }\Phi \left( \lambda ,\varepsilon ,\xi \right) =-\left(
-iB+2\varepsilon \xi \right) \Phi ^{2}\left( \lambda ,\varepsilon ,\xi
\right) .
\end{equation*}

Due to $R$-positivity of the operator $A$ and by assumption (2) the sets 
\begin{equation*}
\left\{ -\left( iB+2\varepsilon \xi \right) \Phi ^{2}\left( \lambda
,\varepsilon ,\xi \right) :\xi \in R\backslash \left\{ 0\right\} \right\}
,\left\{ A\Phi \left( \lambda ,\varepsilon ,\xi \right) :\xi \in R\backslash
\left\{ 0\right\} \right\}
\end{equation*}%
are $R$-bounded. Then in view of the Kahane's contraction principle and from
the product properties of the collection of $R$-bounded operators (see e.g. $%
\left[ \text{4}\right] $ Lemma 3.5, Proposition 3.4) we obtain 
\begin{equation}
\text{ }\sup\limits_{\lambda ,\varepsilon }R\left\{ \xi ^{i}\frac{d}{d\xi
^{i}}\Psi \left( \lambda ,\varepsilon ,\xi \right) \text{: }\xi \in
R\backslash \left\{ 0\right\} \right\} \leq M_{1}\text{, }i=0\text{, }1. 
\tag{2.23}
\end{equation}%
Namely, the $R$-bound of the above sets are independent on $\varepsilon $
and $\lambda $. Next, let us consider $\sigma \left( \lambda ,\varepsilon
,\xi \right) .$ It is clear to see that 
\begin{equation}
\left\Vert \sigma \left( \lambda ,\varepsilon ,\xi \right) \right\Vert
_{B\left( E\right) }\leq C\left\vert \lambda \right\vert
\sum\limits_{j=0}^{2}\left[ \varepsilon ^{\frac{1}{2}}\left\vert \xi
\right\vert \left\vert \lambda \right\vert ^{-\frac{1}{2}}\right]
^{j}\left\Vert \Phi \left( \lambda ,\varepsilon ,\xi \right) \right\Vert
_{B\left( E\right) }.  \tag{2.24}
\end{equation}%
Then by using the well known inequality $y^{j}\leq C\left( 1+y^{m}\right) ,$ 
$y\geq 0,$ $j\leq m$ for $y=\left( \varepsilon ^{\frac{1}{2}}\left\vert
\lambda \right\vert ^{-\frac{1}{2}}\left\vert \xi \right\vert \right) ^{j}$\
and $m=2$ we get the uniform estimate 
\begin{equation}
\left\vert \sum\limits_{j=0}^{2}\varepsilon ^{\frac{j}{2}}\left\vert \lambda
\right\vert ^{1-\frac{j}{2}}\xi ^{j}\right\vert \leq C\left( 1+\varepsilon
\xi ^{2}\left\vert \lambda \right\vert ^{-1}\right) .  \tag{2.25}
\end{equation}%
From $\left( 2.24\right) $ and $\left( 2.25\right) $\ we have the uniform
estimate 
\begin{equation}
\left\Vert \sigma \left( \lambda ,\varepsilon ,\xi \right) \right\Vert
_{B\left( E\right) }\leq C\left\vert \lambda \right\vert \left(
1+\varepsilon \xi ^{2}\left\vert \lambda \right\vert ^{-1}\right) \left(
1+\varepsilon \xi ^{2}+\left\vert \lambda \right\vert \right) ^{-1}\leq C. 
\tag{2.26}
\end{equation}

Due to $R$-positivity of the operator $A,$ the set%
\begin{equation*}
\left\{ \left( \left\vert \lambda \right\vert +\varepsilon \xi ^{2}\right)
\Phi \left( \lambda ,\varepsilon ,\xi \right) :\xi \in R\backslash \left\{
0\right\} \right\}
\end{equation*}%
is $R$-bounded. Then from $\left( 2.26\right) $ and by Kahane's contraction
principle we obtain 
\begin{equation}
\sup\limits_{\lambda ,\varepsilon }R\left\{ \xi ^{i}\frac{d}{d\xi ^{i}}%
\sigma \left( \lambda ,\varepsilon ,\xi \right) :\xi \in R\backslash \left\{
0\right\} \right\} \leq M_{2}\text{, }i=0\text{, }1.  \tag{2.27}
\end{equation}%
By multiplier theorem (see e.g $\left[ 23\right] $) from estimates $\left(
2.23\right) $ and $\left( 2.27\right) $ it follows that $\Psi $ and $\sigma $
are uniform collection of multipliers in $L_{p}\left( R;E\right) .$ Then, by
using the equality\ $\left( 2.21\right) $ we obtain that problem $\left(
2.19\right) $ has a solution $u\in W_{p}^{2}\left( R;E\left( A\right)
,E\right) $ and the uniform estimate holds

\begin{equation}
\sum\limits_{j=0}^{2}\varepsilon ^{\frac{j}{2}}\left\vert \lambda
\right\vert ^{1-\frac{j}{2}}\left\Vert u^{\left( j\right) }\right\Vert
_{L_{p}\left( R;E\right) }+\left\Vert Au\right\Vert _{L_{p}\left( R;E\right)
}\leq C\left\Vert \bar{f}\right\Vert _{L_{p}\left( R;E\right) }.  \tag{2.28}
\end{equation}

Let $u_{1}$ be the restriction of $u$ on $\left( 0,T\right) .$ Then the
estimate $\left( 2.28\right) $ implies that $u_{1}\in Y$. By virtue of
Theorem A$_{3}$ we get 
\begin{equation*}
u_{1}^{\left( m_{k}\right) }\left( .\right) \in \left( E\left( A\right)
;E\right) _{\theta _{k},p},\text{ }k=1,2.
\end{equation*}

Hence, $L_{k}u_{1}\in E_{k}.$ Thus, by Theorem 3.1 problem $\left(
2.20\right) $\ has a unique solution $u_{2}\in Y$ for sufficiently large $%
\left\vert \lambda \right\vert $ and

\begin{equation}
\sum\limits_{j=0}^{2}\varepsilon ^{\frac{j}{2}}\left\vert \lambda
\right\vert ^{1-\frac{j}{2}}\left\Vert u_{2}^{\left( j\right) }\right\Vert
_{X}+\left\Vert Au_{2}\right\Vert _{X}\leq C\sum\limits_{k=1}^{2}\left[
\left\Vert f_{k}\right\Vert _{E_{k}}+\left\vert \lambda \right\vert
^{1-\theta _{k}}\left\Vert f_{k}\right\Vert _{E}+\right.  \notag
\end{equation}

\begin{equation}
\left. \varepsilon ^{\theta _{k}}\left\Vert u_{1}^{\left( m_{k}\right)
}\right\Vert _{C\left( \left[ 0,T\right] ;E_{k}\right) }+\varepsilon
^{\theta _{k}}\left\vert \lambda \right\vert ^{1-\theta _{k}}\left\Vert
u_{1}\right\Vert _{C\left( \left[ 0,T\right] ;E\right) }\right] .  \tag{2.29}
\end{equation}%
Moreover, from $\left( 2.28\right) $ we obtain

\begin{equation}
\sum\limits_{j=0}^{2}\varepsilon ^{\frac{j}{2}}\left\vert \lambda
\right\vert ^{1-\frac{j}{2}}\left\Vert u_{1}^{\left( j\right) }\right\Vert
_{X}+\left\Vert Au_{1}\right\Vert _{X}\leq C\left\Vert f\right\Vert _{X}. 
\tag{2.30}
\end{equation}

Therefore, in virtue of Theorem A$_{3}$ and by estimate $\left( 2.30\right) $
we have

\begin{equation}
\varepsilon ^{\theta _{k}}\left\Vert u_{1}^{\left( m_{k}\right) }\left(
.\right) \right\Vert _{E_{k}}\leq C\left\Vert u_{1}\right\Vert
_{W_{p,\varepsilon }^{2}\left( 0,T;E\left( A\right) ,E\right) }\leq
C\left\Vert f\right\Vert _{L_{p}\left( 0,T;E\right) }.  \tag{2.31}
\end{equation}

In virtue of Theorem A$_{4}$ for $\lambda =\mu ^{2},$ $u\in W_{p}^{2}\left(
0,T;E\right) $ we obtain

\begin{equation}
\ \left\vert \mu \right\vert ^{2-m_{k}}\varepsilon ^{\theta _{k}}\left\Vert
u^{\left( m_{k}\right) }\left( .\right) \right\Vert _{E}\leq C\left[
\left\vert \mu \right\vert ^{\frac{1}{p}}\left\Vert \varepsilon u^{\left(
2\right) }\right\Vert _{X}+\left\vert \mu \right\vert ^{2+\frac{1}{p}%
}\left\Vert u\right\Vert _{X}\right] .  \tag{2.32}
\end{equation}

Hence, from estimates $\left( 2.29\right) $, $\left( 2.31\right) $ and $%
\left( 2.32\right) $ we have

\begin{equation}
\sum\limits_{j=0}^{2}\varepsilon ^{\frac{j}{2}}\left\vert \lambda
\right\vert ^{1-\frac{j}{2}}\left\Vert u_{2}^{\left( j\right) }\right\Vert
_{X}+\left\Vert Au_{2}\right\Vert _{X}\leq  \tag{2.33}
\end{equation}

\begin{equation*}
C\left( \left\Vert f\right\Vert _{X}+\sum\limits_{k=1}^{2}\left( \left\Vert
f_{k}\right\Vert _{E_{k}}+\left\vert \lambda \right\vert ^{1-\theta
_{k}}\left\Vert f_{k}\right\Vert _{E}\right) \right) .
\end{equation*}%
Finally, from $\left( 2.30\right) $ and $\left( 2.33\right) $ we obtain $%
\left( 2.18\right) .$

\begin{center}
\textbf{3. Singular perturbation problem for abstract elliptic equation}
\end{center}

Consider the problem $\left( 1.2\right) $, i.e. the following Cauchy problem
for abstract parabolic equation 
\begin{equation}
Bu^{\prime }\left( t\right) +Au\left( t\right) =f_{0}\left( t\right) ,\text{ 
}t\in \left( 0,T\right) ,  \tag{3.1}
\end{equation}%
\begin{equation}
u\left( 0\right) =u_{0},  \tag{3.2}
\end{equation}%
where $A$, $B$ are linear operators in a Banach space $E.$

The problem $\left( 2.1\right) -\left( 2.2\right) $ can be regarded as the
singular perturbation problem for $\left( 3.1\right) -\left( 3.2\right) .$

In this section we prove the following result:

\textbf{Theorem 3.1. }Let the Condition 2.1 hold and the operator $-AB^{-1}$
generates analytic semigroup in $E.$ Moreover, assume:

( H$_{1}$) $f_{1}\left( \varepsilon \right) \in E,$ $f_{2}\left( \varepsilon
\right) \in D\left( A\right) $, $f_{1}\left( \varepsilon \right) \rightarrow
u_{1}$ in $E$ and $f_{2}\left( \varepsilon \right) \rightarrow u_{0}$ in $%
E\left( A\right) $ as $\varepsilon \rightarrow 0;$

( H$_{2}$) $f\left( \varepsilon ,.\right) \in L_{p}\left( 0,T;E\right) $ and 
$f\left( \varepsilon ,.\right) \rightarrow f_{0}\left( .\right) $ in $X$ as $%
\varepsilon \rightarrow 0.$

Then;

(a) the solution of the \ equation $\left( 2.1\right) $ for $\lambda =0$
converges to the corresponding solution of $(3.1)$ in $X$ as $\varepsilon
\rightarrow 0;$

(b) the solution of $\left( 2.1\right) -\left( 2.2\right) $ converges to the
corresponding solution of $(3.1)-\left( 3.2\right) $ in $E$ as $\varepsilon
\rightarrow 0$ uniformly in $t$ on compact intervals of $\left( 0,T\right) .$

\textbf{Proof. }By virtue of Theorem 2.2, there is a unique solution of $%
\left( 2.1\right) -\left( 2.2\right) $ expressed as%
\begin{equation}
u\left( t,\varepsilon \right) =M\left( t,\varepsilon \right) f_{1}\left(
\varepsilon \right) +N\left( t,\varepsilon \right) f_{2}\left( \varepsilon
\right) +r_{\left[ 0,T\right] }F^{-1}\Phi \left( \xi ,\varepsilon \right) F%
\bar{f}\left( \xi \right) ,  \tag{3.3}
\end{equation}%
where 
\begin{equation*}
M\left( \varepsilon ,t\right) =D^{-1}\left( \varepsilon \right) \left\{
U_{1}\left( t,\varepsilon \right) U_{2}\left( T,\varepsilon \right) \left[
\varepsilon \beta _{1}Q_{2,}\left( \varepsilon \right) +\beta _{0}\right]
\right. -
\end{equation*}%
\begin{equation}
\left. U_{2}\left( t,\varepsilon \right) U_{1}\left( T,\varepsilon \right) 
\left[ \varepsilon \beta _{1}Q_{1}\left( \varepsilon \right) +\beta _{0}%
\right] \right\} ,  \tag{3.4}
\end{equation}%
\begin{equation*}
N\left( \varepsilon ,t\right) =D^{-1}\left( \varepsilon \right) \left\{
U_{2}\left( \varepsilon ,t\right) \left[ \varepsilon \alpha _{1}Q_{1}\left(
\varepsilon \right) +\alpha _{0}\right] -\right.
\end{equation*}%
\begin{equation*}
\left. U_{1}\left( \varepsilon ,t\right) \left[ \varepsilon \alpha
_{1}Q_{2}\left( \varepsilon \right) +\alpha _{0}\right] \right\} ,
\end{equation*}
$\bar{f}$ is a zero extens\i on of $f$ on $\mathbb{R}\setminus \left[ 0,T%
\right] ,$ $r_{\left[ 0,1\right] }$ is a restriction operator from $\mathbb{R%
}$ to $\left[ 0,T\right] $,%
\begin{equation*}
U_{1}\left( x,\varepsilon \right) =\exp -\left\{ xQ_{_{1}}\left( \varepsilon
\right) \right\} ,\text{ }U_{2}\left( x,\varepsilon \right) =\exp -\left\{
xQ_{2}\left( \varepsilon \right) \right\} ,
\end{equation*}

$D^{-1}\left( \varepsilon \right) $, $Q_{_{1}}\left( \varepsilon \right) $, $%
Q_{2}\left( \varepsilon \right) $ are denote $D_{\lambda }^{-1}\left(
\varepsilon \right) ,$ $Q_{_{1,\lambda }}\left( \varepsilon \right) $, $%
Q_{2,\lambda }\left( \varepsilon \right) $ for $\lambda =0$, respectively
and 
\begin{equation*}
\Phi \left( \xi ,\varepsilon \right) =\left( A-i\xi B+\varepsilon \xi
^{2}\right) ^{-1},
\end{equation*}%
Let us show that the solution $u\left( \varepsilon ,.\right) $ of $\left(
2.1\right) -\left( 2.2\right) $ approaches to the corresponding solution of $%
\left( 3.1\right) -\left( 3.2\right) $ in $E$ under conditions (H$_{1}$) and
(H$_{1}$). Since $A$ and $B$ are close operators, it is clear to see that%
\begin{equation*}
\Phi _{0}\left( \xi \right) =\left( A-i\xi B\right) ^{-1}
\end{equation*}%
is a Fourier transform of $Bu^{\prime }\left( t\right) +Au\left( t\right) $
and from $\left( 3.1\right) $ we get that 
\begin{equation*}
\bar{u}\left( t\right) =r_{\left[ 0,T\right] }F^{-1}\Phi _{0}\left( \xi
\right) F\bar{f}_{0}\left( \xi \right) \text{ }
\end{equation*}%
is a solution of the equation $\left( 3.1\right) $, where under Condition
2.1 $\Phi _{0}\left( \xi \right) $\ is uniformly bounded in $\xi \in \mathbb{%
R}.$ The operator functions $\Phi \left( \xi ,\varepsilon \right) $, $\Phi
_{0}\left( \xi \right) $ are uniform bounded and are multipliers in $%
L_{p}\left( R;E\right) $ (see the proof of Theorem 2.2). It is clear to see
that%
\begin{equation}
\Phi \left( \xi ,\varepsilon \right) \rightarrow \Phi _{0}\left( \xi \right) 
\text{ in }B\left( E\right)  \tag{3.5}
\end{equation}%
as $\varepsilon \rightarrow 0$ uniformly in $\xi $ and $\lambda $. Moreover,
we get 
\begin{equation*}
\left\Vert \Phi \left( \xi ,\varepsilon \right) F\bar{f}\left( \xi
,\varepsilon \right) -\Phi _{0}\left( \lambda ,\xi \right) F\bar{f}%
_{0}\left( \xi \right) \right\Vert _{E}\leq
\end{equation*}%
\begin{equation}
\left\Vert \Phi \left( \xi ,\varepsilon \right) F\bar{f}\left( \xi
,\varepsilon \right) -\Phi \left( \xi ,\varepsilon \right) F\bar{f}%
_{0}\left( \xi \right) \right\Vert _{E}+  \tag{3.6}
\end{equation}%
\begin{equation*}
\left\Vert \Phi \left( \xi ,\varepsilon \right) F\bar{f}_{0}\left( \xi
\right) -\Phi _{0}\left( \xi \right) F\bar{f}_{0}\left( \xi \right)
\right\Vert _{E}.
\end{equation*}

Since $\bar{f}\left( \xi ,\varepsilon \right) \rightarrow \bar{f}_{0}\left(
\xi \right) $ in $E$\ as $\varepsilon \rightarrow 0$ for a.e. $\xi \in 
\mathbb{R}$, $\Phi \left( \xi ,\varepsilon \right) $ is bounded in $E$ for
all $\xi \in \mathbb{R}$ and the Fourier transform $F$ is continuous in $X$.
Then we get 
\begin{equation}
\left\Vert \Phi \left( \xi ,\varepsilon \right) F\bar{f}\left( \xi
,\varepsilon \right) -\Phi \left( \xi ,\varepsilon \right) F\bar{f}%
_{0}\left( \xi \right) \right\Vert _{E}\rightarrow 0  \tag{3.7}
\end{equation}%
\ as $\varepsilon \rightarrow 0$ for a.e. for $\xi \in \mathbb{R}$

By the same reason and due to $\Phi \left( \xi ,\varepsilon \right)
\rightarrow \Phi _{0}\left( \xi \right) $ in $B\left( E\right) $\ as $%
\varepsilon \rightarrow 0$ uniformly in $\lambda $ and $\xi ,$ we have

\begin{equation}
\left\Vert \Phi \left( \xi ,\varepsilon \right) F\bar{f}_{0}\left( \xi
\right) -\Phi _{0}\left( \xi \right) F\bar{f}_{0}\left( \xi \right)
\right\Vert _{E}\rightarrow 0.  \tag{3.8}
\end{equation}

Then due to boundedness of $F^{-1}$ from $\left( 3.5\right) -\left(
3.8\right) $ we obtain

\begin{equation*}
\left\Vert F^{-1}\Phi \left( \xi ,\varepsilon \right) F\bar{f}\left( \xi
,\varepsilon \right) -F^{-1}\Phi _{0}\left( \xi \right) F\bar{f}_{0}\left(
\xi \right) \right\Vert _{X}\rightarrow 0
\end{equation*}%
as $\varepsilon \rightarrow 0$, i.e., 
\begin{equation}
r_{\left[ 0,1\right] }F^{-1}\Phi \left( \xi ,\varepsilon \right) F\bar{f}%
\left( \xi ,\varepsilon \right) \rightarrow r_{\left[ 0,1\right] }F^{-1}\Phi
_{0}\left( \xi \right) F\bar{f}_{0}\left( \xi \right) \text{ in }X. 
\tag{3.9}
\end{equation}

We have proved the assertion (a). Now, let us show the assertion (b).
Indeed, known that (see e.g $\left[ \text{1, \S 3 }\right] $, $\left[ \text{%
2, \S\ 1.5}\right] ,$ $\left[ \text{14, \S\ 4.2}\right] $) there is a unique
solution of the Cauchy problem $\left( 3.1\right) -\left( 3.2\right) $ for $%
f\in L_{p}\left( 0,T;E\right) $ expressed as%
\begin{equation*}
u\left( t\right) =U_{0,\lambda }\left( t\right)
u_{0}+\dint\limits_{0}^{t}U_{0,\lambda }\left( t-\tau \right) f_{0}\left(
\tau \right) d\tau ,
\end{equation*}%
where $U_{0,\lambda }\left( t\right) $ is an analytic semigroup in $E$
generated by the operator 
\begin{equation*}
-\text{ }A_{0}\left( \lambda \right) =-A_{\lambda }B^{-1}.
\end{equation*}

Due to uniform boundedness of $D^{-1}\left( \varepsilon \right) $ and by
estimates of analytic semigroups from $\left( 3.4\right) $ we obtain 
\begin{equation*}
\left\Vert M\left( t,\varepsilon \right) f_{1}\right\Vert _{E}\leq C\left\{
\left\Vert U_{2}\left( 1,\varepsilon \right) \right\Vert _{B\left( E\right)
} \left[ \left\Vert U_{1}\left( t,\varepsilon \right) Q\left( \varepsilon
\right) f_{1}\right\Vert _{E}\right. \right. +
\end{equation*}%
\begin{equation*}
\left. \left\Vert U_{1}\left( t,\varepsilon \right) f_{1}\right\Vert _{E} 
\right] +\left\Vert U_{1}\left( 1,\varepsilon \right) \right\Vert _{B\left(
E\right) }\left[ \left\Vert U_{2}\left( t,\varepsilon \right) Q\left(
\varepsilon \right) f_{1}\right\Vert _{E}+\right.
\end{equation*}%
\begin{equation}
\left. \left. \left\Vert U_{2}\left( t,\varepsilon \right) f_{1}\right\Vert
_{E}\right] \right\} \leq C_{1}\exp \left\{ -\varepsilon ^{-1}\omega
t\right\} \left\Vert f_{1}\right\Vert _{E},  \tag{3.12}
\end{equation}%
for $f_{1}\in E$ where, 
\begin{equation*}
Q=Q\left( \varepsilon \right) =\left( B^{2}+4\varepsilon A\right) ^{\frac{1}{%
2}},\text{ }\omega >0.
\end{equation*}

From $\left( 3.4\right) $ in a similar way, for $f_{2}\in E$ we get 
\begin{equation}
\left\Vert N\left( t,\varepsilon \right) f_{2}\right\Vert _{E}\leq C\left\{
\left\Vert U_{1}\left( t,\varepsilon \right) Q\left( \varepsilon \right)
f_{2}\right\Vert _{E}\right. +  \tag{3.13}
\end{equation}%
\begin{equation*}
\left\Vert U_{2}\left( t,\varepsilon \right) Q\left( \varepsilon \right)
f_{2}\right\Vert _{E}+\left. \left\Vert U_{1}\left( t,\varepsilon \right)
f_{2}\right\Vert _{E}\right\} \leq C_{0}\left\Vert f_{2}\right\Vert _{E}.
\end{equation*}

From $\left( 3.12\right) $ and $\left( 3.13\right) $ we have 
\begin{equation}
\lim_{\varepsilon \rightarrow 0}\left\Vert M\left( t,\varepsilon \right)
\right\Vert _{B\left( E\right) }=0,\text{ }\lim_{\varepsilon \rightarrow
0}\left\Vert N\left( t,\varepsilon \right) \right\Vert _{B\left( E\right)
}\leq C_{0}.  \tag{3.14}
\end{equation}

Let us show that 
\begin{equation}
K\left[ N\left( .,\varepsilon \right) -U_{0,}\left( .\right) \right]
\upsilon =U_{0}\ast \left[ N\left( .,\varepsilon \right) -\varepsilon
^{-1}BM\left( .,\varepsilon \right) \right] A_{0}\upsilon  \tag{3.15}
\end{equation}%
for all $\upsilon \in D\left( A_{0}\right) ,$ where $K$ is a uniform bounded
operator in $E.$

Indeed, the Laplace transform of $U_{0}\left( .\right) $, $U_{1}\left(
.,\varepsilon \right) $, $U_{2}\left( .,\varepsilon \right) $ gives the
resolvent $R(s,A_{0})$, $R(s,B+Q)$, $R(s,B+Q)$, respectively. Hence, by
using the linearity and convolution properties of the Laplace transform, $%
\left( 3.15\right) ,$ $\left( 3.4\right) $ and $\left( 2.6\right) $ it
sufficient to show

\begin{equation*}
KD^{-1}\left[ \left( \varepsilon \alpha _{1}Q_{1}+\alpha _{0}\right) R\left(
s,Q_{2}\right) -\left( \varepsilon \alpha _{1}Q_{2}+\alpha _{0}\right)
R\left( s,Q_{1}\right) \right] -K_{\lambda }R\left( s,A_{0}\right) =
\end{equation*}%
\begin{equation}
A_{0}R\left( s,A_{0}\right) \left\{ D_{\lambda }^{-1}\left[ \left(
\varepsilon \alpha _{1}Q_{1}+\alpha _{0}\right) \right] R\left(
s,Q_{2}\right) \right. -  \tag{3.16}
\end{equation}%
\begin{equation*}
\left. \left( \varepsilon \alpha _{1}Q_{2}+\alpha _{0}\right) R\left(
s,Q_{1}\right) \right] -\varepsilon ^{-1}BD^{-1}\left[ \left( \varepsilon
\beta _{1}Q_{2}+\beta _{0}\right) U_{2}\left( \varepsilon ,T\right) R\left(
s,Q_{1}\right) \right. +
\end{equation*}%
\begin{equation*}
\left. \left. \left( \varepsilon \beta _{1}Q_{1}+\beta _{0}\right)
U_{1}\left( \varepsilon ,T\right) R\left( s,Q_{2}\right) \right] \right\} .
\end{equation*}

Indeed, by using $\left( 2.6\right) $, the resolvent equation, the
exponential properties of strongly continuous semigroups we get that there
is a bounded operator $K$ in $E$ that $\left( 3.16\right) $ is satisfied.$.$
Hence, from $\left( 3.4\right) $ and $\left( 3.13\right) $ for $\upsilon \in
D\left( A\right) $ we get%
\begin{equation*}
\left\Vert \left[ N\left( .,\varepsilon \right) -U_{0}\left( .\right) \right]
\upsilon \right\Vert _{E}\leq C_{1}\exp \left\{ -\varepsilon ^{-1}\omega
t\right\} \left\Vert A_{0}\upsilon \right\Vert _{E}+
\end{equation*}%
\begin{equation}
C_{2}\exp \left\{ -\varepsilon ^{-1}\omega t\right\} \left\Vert U_{0}\left(
.\right) \right\Vert _{B\left( E\right) }\left\Vert A_{0}\upsilon
\right\Vert _{E}.  \tag{3.17}
\end{equation}

Then from $\left( 3.3\right) $, $\left( 3.4\right) $ and $\left( 3.17\right) 
$ for $f_{1}\in E$, $f_{2}\in D\left( A\right) $ we deduced 
\begin{equation*}
\left\Vert u\left( .,\varepsilon \right) -u\left( .\right) \right\Vert
_{E}\leq \left\Vert M\left( .,\varepsilon \right) f_{1}\right\Vert
_{E}+\left\Vert N\left( .,\varepsilon \right) f_{2}-U_{0}\left( .\right)
u_{0}\right\Vert _{E}+
\end{equation*}%
\begin{equation}
\left\Vert f\left( .,\varepsilon \right) -f_{0}\left( .\right) \right\Vert
_{E}\leq C_{1}\exp \left\{ -\varepsilon ^{-1}\omega t\right\} \left\Vert
f_{1}\right\Vert _{E}+  \tag{3.18}
\end{equation}%
\begin{equation*}
C_{2}\exp \left\{ -\varepsilon ^{-1}\omega t\right\} \left\Vert
f_{2}\right\Vert _{E}+\left\Vert f\left( \varepsilon ,.\right) -f_{0}\left(
.\right) \right\Vert _{E},\text{ }
\end{equation*}%
\ 

By conditions (H$_{1}$) and (H$_{2}$) we get 
\begin{equation*}
\exp \left\{ -\varepsilon ^{-1}\omega t\right\} \rightarrow 0\text{ as\ \ }%
\varepsilon \rightarrow 0
\end{equation*}%
uniformly with respect to $t$ on all compact $\sigma \subset \left(
0,T\right) $. Then from $\left( 3.18\right) $ we obtain the assertion.

\begin{center}
\textbf{4.} \textbf{Wentzell-Robin type mixed problem for elliptic equation}
\end{center}

Consider the BVP $\left( 1.4\right) -\left( 1.5\right) .$ For $\mathbf{p=}%
\left( p,2\right) $ and $L_{\mathbf{p}}\left( \Omega \right) $ will denote
the space of all $\mathbf{p}$-summable scalar-valued\ functions with mixed
norm. Analogously, $W_{\mathbf{p}}^{2}\left( \Omega \right) $ denotes the
Sobolev space with corresponding mixed norm, i.e., $W_{\mathbf{p}}^{2}\left(
\Omega \right) $ denotes the space of all functions $u\in L_{\mathbf{p}%
}\left( \Omega \right) $ possessing the derivatives $\frac{\partial ^{2}u}{%
\partial x^{2}},$ $\frac{\partial ^{2}u}{\partial y^{2}}\in L_{\mathbf{p}%
}\left( \Omega \right) $ with the norm 
\begin{equation*}
\ \left\Vert u\right\Vert _{W_{\mathbf{p}}^{2}\left( \Omega \right)
}=\left\Vert u\right\Vert _{L_{\mathbf{p}}\left( \Omega \right) }+\left\Vert 
\frac{\partial ^{2}u}{\partial x^{2}}\right\Vert _{L_{\mathbf{p}}\left(
\Omega \right) }+\left\Vert \frac{\partial ^{2}u}{\partial y^{2}}\right\Vert
_{L_{\mathbf{p}}\left( \Omega \right) }.
\end{equation*}

\textbf{Condition 4.1 }Assume:

(1) $K\left( .,.\right) \in C\left( \left[ 0,T\right] \times \left[ 0,1%
\right] \right) ;$

(2)\ $a$ is positive, $b$ is a real-valued functions on $\left( 0,1\right) ;$

(3) $a\left( .\right) \in C\left( 0,1\right) $ and%
\begin{equation*}
\exp \left( -\dint\limits_{\frac{1}{2}}^{x}b\left( t\right) a^{-1}\left(
t\right) dt\right) \in L_{1}\left( 0,1\right) .
\end{equation*}

\bigskip\ In this section, we present the following result:

\bigskip \bigskip \textbf{Theorem 4.1. }Suppose the Condition 4.1 hold. Then:

(a) For $f\in L_{\mathbf{p}}\left( \Omega \right) $, $p$, $p_{1}\in \left(
1,\infty \right) $ problem $\left( 1.4\right) -\left( 1.5\right) $ has a
unique solution $u\in W_{\mathbf{p}}^{2}\left( \Omega \right) $\ and the
following uniform coercive estimate holds 
\begin{equation*}
\sum\limits_{i=0}^{2}\varepsilon ^{\frac{i}{2}}\left\vert \lambda
\right\vert ^{1-\frac{i}{2}}\left\Vert \frac{\partial ^{i}u}{\partial x^{i}}%
\right\Vert _{L_{\mathbf{p}}\left( \Omega \right) }+\left\Vert \frac{%
\partial ^{2}u}{\partial y^{2}}\right\Vert _{L_{\mathbf{p}}\left( \Omega
\right) }+\ \left\Vert u\right\Vert _{L_{\mathbf{p}}\left( \Omega \right)
}\leq
\end{equation*}

\begin{equation*}
C\left[ \left\Vert f\right\Vert _{L_{\mathbf{p}}\left( \Omega \right)
}+\left\Vert f_{1}\right\Vert _{L_{p}\left( 0,1\right) }+\left\Vert
f_{2}\right\Vert _{W_{p}^{1}\left( 0,1\right) }\right] ;
\end{equation*}%
(b) the solution of the equation $\left( 1.4\right) $ for $\lambda =0$
converges to the corresponding solution of the following equation%
\begin{equation*}
-\left( a\frac{\partial ^{2}u}{\partial y^{2}}+b\frac{\partial u}{\partial y}%
\right) +\dint\limits_{0}^{1}K\left( y,\tau \right) \frac{\partial }{%
\partial t}u\left( t,y,\tau \right) d\tau =f\left( t,y\right) \text{, }
\end{equation*}%
in $L_{\mathbf{p}}\left( \Omega \right) $ as $\varepsilon \rightarrow 0;$

(c) the solution of $\left( 1.4\right) -\left( 1.5\right) $ converges to the
corresponding solution of the following mixed problem 
\begin{equation*}
-\left( a\frac{\partial ^{2}u}{\partial y^{2}}+b\frac{\partial u}{\partial y}%
\right) +\dint\limits_{0}^{1}K\left( y,\tau \right) \frac{\partial }{%
\partial t}u\left( t,y,\tau \right) d\tau =f\left( t,y\right) \text{, }
\end{equation*}%
\begin{equation*}
u\left( 0,y\right) =0\text{ for a.e. }y\in \left( 0,1\right) ,
\end{equation*}%
\ \ \ 
\begin{equation*}
a\left( j\right) u_{yy}\left( t,j,\varepsilon \right) +b\left( j\right)
u_{y}\left( t,j,\varepsilon \right) =0\text{, }j=0,1\text{ for a.e. }t\in
\left( 0,T\right) ,
\end{equation*}%
in $L_{p}\left( 0,1\right) $ as $\varepsilon \rightarrow 0$ uniformly in $t$
on compact intervals of $\left( 0,T\right) .$

\ \textbf{Proof.} Let $E=L_{2}\left( 0,1\right) $. It is known $\left[ 5%
\right] $\ that $L_{2}\left( 0,1\right) $ is an $UMD$ space. Consider the
operator $A$ defined by 
\begin{equation*}
D\left( A\right) =W_{2}^{2}\left( 0,1;A\left( j\right) u=0\text{, }%
j=0,1\right) ,\text{ }Au=-a\frac{\partial ^{2}u}{\partial y^{2}}+b\frac{%
\partial u}{\partial y}.
\end{equation*}

Therefore, the problem $\left( 1.4\right) -\left( 1.5\right) $ can be
rewritten in the form of $\left( 2.2\right) $, where $u\left( t\right)
=u\left( t,.\right) ,$ $f\left( t\right) =f\left( t,.\right) $\ are
functions with values in $E=L_{2}\left( 0,1\right) .$ By virtue of $\left[ 
\text{8}\right] $ the operator $A$ generates analytic semigroup in $%
L_{2}\left( 0,1\right) $. Then in view of Hill-Yosida theorem (see e.g. $%
\left[ \text{22, \S\ 1.13}\right] $) this operator is sectorial in $%
L_{2}\left( 0,1\right) .$ Since all uniform bounded set in Hilbert apace is
an $R$-bounded (see $\left[ 3\right] $ ), i.e. we get that the operator $A$
is $R$-sectorial in $L_{2}\left( 0,1\right) .$ Then from Theorem 2.2 and
Theorem 3.1 we obtain the assertion.

\begin{center}
\textbf{Acknowledgements}
\end{center}

The author is thanking the library manager of Okan University Kenan \"{O}%
ztop for his help in finding the necessary articles and books in my research
area.

\begin{center}
\ \textbf{References}
\end{center}

\begin{quote}
\ \ \ \ \ \ \ \ \ \ \ \ \ \ \ \ \ \ \ \ \ \ \ \ \ \ \ \ \ \ \ \ \ \ \ \ \ \
\ \ \ \ \ \ \ \ \ \ \ \ \ \ \ \ \ \ \ \ \ \ \ \ \ \ \ \ \ \ \ \ \ \ \ \ \ \
\ \ \ \ \ \ \ \ \ \ 
\end{quote}

\begin{enumerate}
\item Amann H., Linear and quasi-linear equations,1, Birkhauser, Basel 1995.

\item Ashyralyev A., Cuevas C., and Piskarev S., On well-posedness of
difference schemes for abstract elliptic problems in spaces, Numer. Func.
Anal. Opt., (2008)29, (1-2), 43-65.

\item Burkholder D. L., A geometrical conditions that implies the existence
certain singular integral of Banach space-valued functions, Proc. conf.
Harmonic analysis in honor of Antonu Zigmund, Chicago, 1981,Wads Worth,
Belmont, (1983), 270-286.

\item Dore C. and Yakubov S., Semigroup estimates and non coercive boundary
value problems, Semigroup Forum 60 (2000), 93-121.

\item Denk R., Hieber M., Pr\"{u}ss J., $R$-boundedness, Fourier multipliers
and problems of elliptic and parabolic type, Mem. Amer. Math. Soc. 166
(2003), n.788.

\item Engel, K.-J., On singular perturbations on second order Cauchy
problems. Pacific J. Math. 152 (1) (1992), 79-91.

\item Fattorini, H, O, \ The Cauchy Problems, Addison-Wesley, Reading,
Mass., 1983.

\item A. Favini, G. R. Goldstein, J. A. Goldstein and S. Romanelli,
Degenerate second order differential operators generating analytic
semigroups in $L_{p}$ and $W^{1,p}$, Math. Nachr. 238 (2002), 78 --102.

\item Goldstain J. A., Semigroups of linear operators and applications,
Oxfard University Press, Oxfard 1985.

\item Lunardi A., Analytic semigroups and optimal regularity in parabolic
problems, Birkhauser, Verlag, Basel,1995.

\item Krein S. G., Linear differential equations in Banach space, American
Mathematical Society, Providence, 1971.

\item Kisy\`{n}ski J., Sur les \'{e}quations hyperboliques avec petite param%
\'{e}tre. Colloq. Math. 10 (1963), 331-343.

\item Lions J. L and Peetre J., Sur une classe d'espaces d'interpolation,
Inst. Hautes Etudes Sci. Publ. Math., 19(1964), 5-68.

\item Pazy A., Semigroups of Linear Operators and Applications to Partial
Differential Equations. Springer, Berlin, 1983

\item Ragusa M. A., Embeddings for Lorentz-Morrey spaces, J. Optim. Theory
Appl., 154(2)(2012), 491-499.

\item Shahmurov R., Solution of the Dirichlet and Neumann problems for a
modified Helmholtz equation in Besov spaces on an annuals, J. Differential
Equations, 249(3) (2010), 526-550.

\item Shakhmurov V. B., Linear and nonlinear abstract equations with
parameters, Nonlinear Anal., 73(2010), 2383-2397.

\item Shakhmurov V. B., Estimates of approximation numbers for embedding
operators and applications, Acta. Math. Sin., (Engl. Ser.), (2012), 28 (9),
1883-1896.

\item Shakhmurov V. B., Coercive boundary value problems for regular
degenerate differential-operator equations, J. Math. Anal. Appl., 292 (2)
(2004), 605-620.

\item Shakhmurov V. B., Embedding theorems and\ maximal regular differential
operator equations in Banach-valued function spaces, J. Inequal. Appl.,
2(4)(2005), 329-345.

\item Shahmurov R., On strong solutions of a Robin problem modeling heat
conduction in materials with corroded boundary, Nonlinear Anal. Real World
Appl., 13(1) (2011), 441-451.

\item Triebel H., Interpolation theory, Function spaces, Differential
operators, North-Holland, Amsterdam, 1978.

\item Weis L, Operator-valued Fourier multiplier theorems and maximal $L_{p}$
regularity, Math. Ann. 319, (2001), 735-758.

\item Yakubov S. and Yakubov Ya., Differential-operator Equations. Ordinary
and Partial Differential Equations, Chapman and Hall /CRC, Boca Raton, 2000.
\end{enumerate}

\begin{quote}
\ 

\bigskip
\end{quote}

\end{document}